\documentclass[12pt]{amsart}
\usepackage{amscd,amssymb,amsthm,verbatim,tikz-cd,graphicx}
\graphicspath{ {./images/} }
\newcommand{\C}{{\mathbb C}}

\newcommand{\ch}{\operatorname{ch}}

\newcommand{\Coker}{\operatorname{Coker}}

\newcommand{\const}{\operatorname{const.}}

\newcommand{\diam}{\operatorname{diam}}

\newcommand{\Dom}{\operatorname{Dom}}
\newcommand{\dvol}{\operatorname{dvol}}

\newcommand{\HH}{\operatorname{H}}

\newcommand{\Id}{\operatorname{Id}}

\newcommand{\Image}{\operatorname{Im}}

\newcommand{\Int}{\operatorname{int}}

\newcommand{\KK}{\operatorname{K}}

\newcommand{\Ker}{\operatorname{Ker}}

\newcommand{\pt}{\operatorname{pt}}
\newcommand{\Q}{{\mathbb Q}}
\newcommand{\R}{{\mathbb R}}

\newcommand{\Ric}{\operatorname{Ric}}
\newcommand{\Riem}{\operatorname{Riem}}

\newcommand{\Rm}{\operatorname{Rm}}

\newcommand{\tr}{\operatorname{tr}}

\newcommand{\vol}{\operatorname{vol}}
\newcommand{\Z}{{\mathbb Z}}

\numberwithin{equation}{section}
\theoremstyle{plain}

\newtheorem{definition}[equation]{Definition}
\newtheorem{assumption}[equation]{Assumption}
\newtheorem{lemma}[equation]{Lemma}
\newtheorem{theorem}[equation]{Theorem}
\newtheorem{proposition}[equation]{Proposition}
\newtheorem{corollary}[equation]{Corollary}

\newtheorem{conjecture}[equation]{Conjecture}
\newtheorem{hypothesis}[equation]{Hypothesis}

\theoremstyle{remark}

\newtheorem{remark}[equation]{Remark}
\newtheorem{example}[equation]{Example}
\usepackage{graphicx}
\input{amssym.def}
\input{amssym}
\input{epsf}
\usepackage{a4wide}

\begin{document}

\title[Some obstructions to positive scalar curvature on a noncompact manifold]
      {Some obstructions to positive scalar curvature on a noncompact manifold}

\author{John Lott}
\address{Department of Mathematics\\
University of California, Berkeley\\
Berkeley, CA  94720-3840\\
USA
}

\date{September 18, 2025}

\email{lott@berkeley.edu}

\begin{abstract}
We give obstructions for a noncompact manifold to admit a complete Riemannian metric with (nonuniformly)
positive scalar curvature.  We treat both the finite volume and infinite volume cases.
\end{abstract}

\maketitle

\noindent

\section{Introduction} \label{sect1}

There are many results on whether a given compact manifold admits a Riemannian metric with positive scalar curvature (psc).
In this paper we focus on the noncompact case.

An open conjecture for compact manifolds says that an aspherical compact smooth
manifold cannot admit a psc
metric.
There are two main approaches to this conjecture. The first one uses minimal hypersurfaces,
following the work of Schoen-Yau \cite{Schoen-Yau (1979)},
and $\mu$-bubbles as introduced by Gromov \cite[Section $5 \frac56$]{Gromov (1996)}.  Recent advances are by
Chodosh-Li \cite{Chodosh-Li (2024)} and Gromov \cite{Gromov (2020)}.  The other approach, which we follow, uses Dirac operators.

The above conjecture has an extension to compact manifolds that may not be aspherical.  If $M$
is a compact connected oriented $n$-dimensional smooth manifold, choose a basepoint $m_0$ and consider the fundamental group
$\Gamma = \pi_1(M, m_0)$. There is a pointed connected CW-complex $B\Gamma$ with the property that 
$\pi_1(B\Gamma) = \Gamma$ and the universal cover of $B\Gamma$ is contractible. There is a
classifying map $\nu : M \rightarrow B\Gamma$, unique up to homotopy, that induces an isomorphism
on $\pi_1$.  If $[M] \in 
\HH_n(M; \Q)$ is the fundamental class in rational homology then the extended conjecture says that
nonvanishing of the pushforward
$\nu_*[M] \in \HH_n(B\Gamma; \Q)$ is an obstruction for $M$ to admit a psc metric.
If $M$ is aspherical then one recovers the previous conjecture.  There are many results on this
extended conjecture, using Dirac operators \cite{Rosenberg (2007)}.

This paper is concerned with obstructions to complete psc metrics on noncompact manifolds.  
There is also a long history to this problem, going back to the Gromov-Lawson paper \cite{Gromov-Lawson (1983)}.
There is a technical difference between looking at metrics with uniformly positive scalar
curvature, i.e. bounded below by a positive constant, and metrics with nonuniformly positive scalar curvature.  In the first case the
Dirac operator is Fredholm, while in the second case it need not be Fredholm.  We are
interested in psc metrics that may not have uniformly positive scalar curvature; the study
of such metrics goes back to \cite[Section 6]{Gromov-Lawson (1983)}.
Our motivation comes from conjectures relating scalar curvature to simplicial volume, for
compact manifolds.

Here is a test question.  Suppose that $Y$ is a compact connected oriented $n$-dimensional
manifold-with-boundary, with
connected boundary $\partial Y$. Choosing a basepoint $y_0 \in \partial Y$, put
$\Gamma = \pi_1(Y, y_0)$ and $\Gamma^\prime = \pi_1(\partial Y, y_0)$.  There is a classifying
map of pairs $\nu : (Y, \partial Y) \rightarrow (B\Gamma, B\Gamma^\prime)$, unique up to
homotopy.  Let $[Y, \partial Y] \in \HH_n(Y, \partial Y; \Q)$ be the fundamental class.
Is nonvanishing of the pushforward $\nu_* [Y, \partial Y] \in \HH_n(B\Gamma, B\Gamma^\prime; \Q)$ an obstruction for the 
interior $\Int(Y) = Y - \partial Y$ of $Y$
to admit a complete psc metric, provided that\\
(a) The homomorphism $\Gamma^\prime \rightarrow \Gamma$ is injective, or \\
(b) The metric has finite volume?

Here if $\Gamma^\prime \rightarrow \Gamma$ is not injective then we define $\HH_n(B\Gamma, B\Gamma^\prime; \Q)$
using the algebraic mapping cone complex; it could more accurately be written as
$\HH_n(B\Gamma^\prime \rightarrow B\Gamma; \Q)$. 

One needs some condition like (a) or (b),
as can be seen if $Y = D^2$. Then $\pi_1(\partial Y) \rightarrow \pi_1(Y)$ is not injective,
$\nu_* [Y, \partial Y] \neq 0$ and
$\Int(Y)$ does admit a complete psc metric, such as a paraboloid, but not one of finite volume.

A special case of the above question is when $Y$ and $\partial Y$ are aspherical, in which 
case $\nu_* [Y, \partial Y]$ is automatically nonzero.

Regarding (a), if the homomorphism $\Gamma^\prime \rightarrow \Gamma$ is not injective then
the obstruction should lie in $\HH_n(B\Gamma, B\widetilde{\Gamma}^\prime; \Q)$, where
$\widetilde{\Gamma}^\prime$ is the image of $\Gamma^\prime$ in $\Gamma$.

We give results in the direction of (a) and (b).  One
main tool is almost flat bundles in the relative setting.  Almost flat bundles were
introduced by Connes-Gromov-Moscovici \cite{Connes-Gromov-Moscovici (1990)} and give obstructions for compact spin
manifolds to have psc metrics.  We review this material in Section \ref{sect2}.  Almost flat
bundles in the relative setting were introduced by Kubota \cite{Kubota (2022)}.  There are actually two
versions: almost flat relative bundles and almost flat stable relative bundles.  They
are relevant for (a) and (b), respectively. 

Our other main technical tool is Callias-type Dirac
operators, as were used for example by Cecchini-Zeidler in 
\cite{Cecchini-Zeidler (2021),Cecchini-Zeidler (2021b)}. This allows us to give localized
obstructions to positive scalar curvature, that apply to incomplete manifolds.
Some of the statements involve the mean curvature of a boundary.

For notation, $R$ denotes scalar curvature.
Our convention for mean curvature is such that $S^{n-1} \subset {D^n}$ has mean curvature
$H = n-1$.

The geometric setup for our localized statements is that we have a region in a manifold
that is assumed to have positive scalar curvature, then there is
an annulus around it with quantitatively positive scalar curvature, followed by a larger annulus
of a certain size that can have some
negative scalar curvature.  What happens outside of the second annulus doesn't matter. 
The precise assumption is the following.
(See Figure 1.)

 \begin{assumption} \label{1.1}
    Given $r_0, D > 0$, put $r_0^\prime = \frac{1}{256} r_0^2 D^2$
and $D^\prime = D + \frac{32}{r_0 D}$. 
Let $M$ be a connected Riemannian spin manifold, possibly with boundary and possibly incomplete. Let $K$ be a compact subset of $M$
containing $\partial M$.
Suppose that 
\begin{itemize}
\item The distance neighborhood $N_{D^\prime}(K)$ lies in a compact submanifold-with-boundary ${\mathcal C}$,
\item $R > 0$ on $K$,
\item $R \ge r_0$ on $N_{D}(K) - K$ and
\item $R \ge - r_0^\prime$ on $N_{D^\prime}(K) - N_{D}(K)$.
\end{itemize}
    \end{assumption}

    \includegraphics[scale=1]{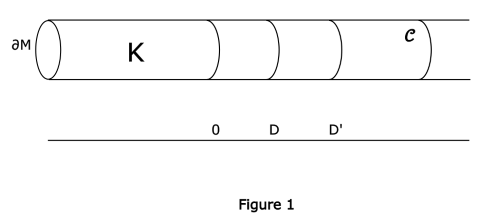}

The first main result just uses almost flat $\KK$-theory classes, denoted by $\KK^*_{af}(\cdot)$. 
    
\begin{theorem} \label{1.2}
Suppose that Assumption \ref{1.1} holds, where
$\partial M$ has nonnegative mean curvature.
Given $\beta \in 
\KK^{-1}_{af}({\mathcal C})$ if $M$ is even dimensional, or
$\beta \in 
\KK^{0}_{af}({\mathcal C})$ if $M$ is odd dimensional, we have
\begin{equation} \label{1.3}
\int_{\partial M} \widehat{A}(T\partial M) \wedge \ch \left(
\beta \big|_{\partial M} \right) = 0.
\end{equation}
\end{theorem}

Our main application of Theorem \ref{1.2} is to finite volume complete manifolds $M$ of dimension at most seven.
Proposition \ref{2.41} says that
there is an exhaustion of $M$ by compact submanifolds $K_i$
whose boundaries $\partial K_i$ have nonnegative mean curvature as seen from ${M - K_i}$. 
(It would be interesting if the dimension restriction could be removed.)
Rather than applying Theorem \ref{1.2} to $K_i$, we apply it to a suitable compact region of $\overline{M - K_i}$
containing $\partial K_i$. 
In this way we obtain end obstructions to the existence of finite volume psc metrics.
As a simple example, there is no complete finite volume psc metric on $[0, \infty) \times T^{n-1}$ if $n \le 7$.

The next main result uses almost flat relative $\KK$-theory classes, denoted by $\KK^*_{af}(\cdot, \cdot)$.
The geometric setup is similar to the previous one, except that now there is no boundary.

\begin{theorem} \label{1.4}
Suppose that Assumption \ref{1.1} holds, where $\partial M = \emptyset$ and
$K$ is a compact submanifold-with-boundary in $M$.
Then given $\beta \in \KK^*_{af} 
\left({\mathcal C},  {\mathcal C} - \Int(K)
\right)$,
we have
\begin{equation} \label{1.5}
\int_{\mathcal C} \widehat{A}(TM) \wedge 
\ch (\beta) = 0.
\end{equation}
\end{theorem}

Theorem \ref{1.4} is relevant to part (a) of the test question above; see Corollary \ref{3.41}.

The third main result combines Theorems \ref{1.2} and \ref{1.4}.
It uses almost flat stable relative $\KK$-theory classes, denoted by $\KK^*_{af,st}(\cdot, \cdot)$
The generators of $\KK^0_{af,st}(\cdot, \cdot)$ differ from generators of
$\KK^0_{af}(\cdot, \cdot)$ essentially by having additional $\KK^{-1}$-generators for the second factor.
This meshes well with the compact exhaustions of finite volume manifolds.
In the geometric assumptions,
it is now assumed that the boundary of the inner compact region has nonnegative mean
curvature as seen from the complement.

\begin{theorem} \label{1.6}
Suppose that Assumption \ref{1.1} holds, where $\partial M = \emptyset$,
$K$ is a compact codimension-zero submanifold-with-boundary in $M$, and $\partial K$ has 
nonnegative mean curvature as seen from
$M-K$.
Then given $\beta \in \KK^{*}_{af,st} \left( 
{\mathcal C},  {\mathcal C} - \Int(K) \right)$,
we have
\begin{equation} \label{1.7}
\int_{\mathcal C} \widehat{A}(TM) \wedge \ch(\beta) = 0.
\end{equation}
\end{theorem}

When combined with the result about compact exhaustions of complete finite volume Riemannian manifolds, 
Theorem \ref{1.6} is relevant to part (b) of the test question above; see Corollary \ref{4.15}.

For technical reasons, our results use various types of almost flat bundles.  Another approach toward 
index theoretic obstructions to positive scalar curvature uses group $C^*$-algebras.  It would be interesting
if similar results could be proved using relative $C^*$-algebras and Mishchenko bundles.

Our interest in psc metrics comes from a conjecture about the
simplicial volume of manifolds having almost nonnegative scalar curvature, with respect to a normalized volume.
We wanted to see if the conjecture can be verified under additional geometric bounds.  We were only
partly successful in this; see Appendix A.

The structure of the paper is the following.  In Section \ref{sect2} we prove Theorem \ref{1.2}
and give an application to finite volume Riemannian manifolds.
In Section \ref{sect3} we prove Theorem \ref{1.4}. We also give obstructions for a compact
manifold-with-boundary to have a psc metric with nonnegative mean curvature on the
boundary. In Section \ref{sect4} we prove Theorem \ref{1.6}.
The appendix has a discussion of simplicial volume. More detailed descriptions are at the
beginnings of the sections.

I thank Clara L\"oh and Antoine Song for helpful discussions, and the referee for useful comments.

\section{Almost flat bundles and end obstructions} \label{sect2}

In this section we use almost flat bundles to construct obstructions for noncompact manifolds to admit complete finite volume
psc metrics.  In Subsection \ref{subsect2.1} we review almost flat bundles and their use for compact manifolds.
Subsection \ref{subsect2.2} gives an obstruction for a manifold-with-boundary 
to have a metric with positive scalar curvature in a neighborhood of the boundary,
with nonnegative mean curvature on the boundary. In Subsection \ref{subsect2.3} we give obstructions for a
manifold to have a complete finite volume psc metric.

\subsection{Almost flat bundles} \label{subsect2.1}
If $X$ is a compact topological space then
generators of $\KK^0(X)$ are 
$\Z_2$-graded vector bundles on $X$, and there is a Chern character
$\ch : \KK^0(X) \rightarrow \HH^{even}(X; \Q)$.
Generators of $\KK^{-1}(X)$ are given
by pairs $(V, \sigma)$ where $V$ is a 
vector bundle on $X$ and $\sigma$ is an
automorphism of $V$. 
There is a Chern character
$\ch : \KK^{-1}(X) \rightarrow \HH^{odd}(X; \Q)$.

If $X$ is noncompact then we take $\KK(X)$ to be the representable K-group, i.e. $\KK^0(X) = [X, \Z \times BU]$.
In what follows, all manifolds will be taken
to be smooth. 

\begin{definition} \label{2.1} \cite{Connes-Gromov-Moscovici (1990)}
Let ${\mathcal N}$ be a compact manifold-with-boundary.  Put a Riemannian metric on ${\mathcal N}$.
Given $\beta \in \KK^0({\mathcal N})$ we say that $\beta$ is almost flat
if for every $\epsilon > 0$, 
we can find 
\begin{itemize}
\item A $\Z_2$-graded Hermitian vector bundle $V$ on ${\mathcal N}$ representing $\beta$, and
\item A Hermitian connection $\nabla^\pm$ on $V$ whose curvature satisfies $\| F^\pm \| \le \epsilon$.
\end{itemize}
\end{definition}

Whether $\beta$ is almost flat is independent of the choice of Riemannian metric on ${\mathcal N}$.
Let $\KK^0_{af}({\mathcal N})$ denote the almost flat elements of $\KK^0({\mathcal N})$. The relevance of almost
flat bundles for questions of positive scalar curvature
comes from the following result.
\begin{theorem} \label{2.2} \cite{Connes-Gromov-Moscovici (1990)}
If $M$ is a closed even dimensional
spin manifold with a 
positive scalar curvature metric, and $\beta \in
\KK^0_{af}(M)$, then
$\int_M \widehat{A}(TM) \wedge
\ch(\beta) = 0$.
\end{theorem}

Recall that a discrete group $\Gamma$ has a classifying space $B\Gamma$ which can be
taken to be a connected CW-complex,
defined up to homotopy, with the property that $\pi_1(B\Gamma) \cong \Gamma$ and the universal cover of $B\Gamma$ is
contractible.
If ${\mathcal N}$ is
connected, let $\nu_{\mathcal N} : {\mathcal N} \rightarrow B\pi_1({\mathcal N}, n_0)$ be the classifying map
for the universal cover of ${\mathcal N}$. It
is known that elements of $\KK^0_{af}({\mathcal N})$
pull back from $B\pi_1({\mathcal N}, n_0)$, in an appropriate sense
\cite{Hunger (2019),Mishchenko-Teleman (2005)}.
This
motivates the following definition.

\begin{definition} \label{2.3}
Let $\Gamma$ be a discrete group.  Given $\eta \in \KK^0(B\Gamma)$, we say that $\eta$ is almost flat if for every 
compact manifold-with-boundary ${\mathcal N}$ and 
every continuous map $\nu : {\mathcal N} \rightarrow B\Gamma$, the pullback $\nu^* \eta$ 
is almost flat.
\end{definition}

Let $\KK^0_{af}(B\Gamma)$ denote the almost flat classes. It is a subgroup of 
$\KK^0(B\Gamma)$.
One can give a
more intrinsic description of $\KK^0_{af}(B\Gamma)$ in terms of
$\Gamma$ \cite{Hunger (2019)}. It seems conceivable that for every discrete group $\Gamma$, the
map $\KK^0_{af}(B\Gamma) \otimes \Q \rightarrow \KK^0(B\Gamma) \otimes \Q$ is an isomorphism,
at least if $B\Gamma$ can be represented by a finite CW-complex. (Our definitions avoid some subtleties in the case of infinite complexes.)
Some examples of groups $\Gamma$ where this is known are listed in \cite{Connes-Gromov-Moscovici (1990)}.

Note that if $B\Gamma$ can be represented by a smooth compact manifold then $\eta \in \KK^0(B\Gamma)$ is almost flat in the sense of Definition \ref{2.3} if and only if it is almost flat
in the sense of Definition \ref{2.1}.

We now give the odd degree version of the previous definitions.

\begin{definition} \label{2.4}
Let ${\mathcal N}$ be a compact manifold-with-boundary. Put a Riemannian metric on ${\mathcal N}$.
Given $\beta \in \KK^{-1}({\mathcal N})$, we say that $\beta$ is almost flat if for every 
$\epsilon > 0$, 
we can find 
\begin{itemize}
\item A Hermitian vector bundle $V$ on ${\mathcal N}$ and a unitary automorphism
$\sigma$ of $V$ so that $(V, \sigma)$ represents $\beta$, and
\item A Hermitian connection $\nabla$ on $V$ whose curvature satisfies $\| F \| \le \epsilon$, and also $\| \nabla \sigma \| \le \epsilon$.
\end{itemize}
\end{definition}

Whether $\beta$ is almost flat is independent of the choice of Riemannian metric on ${\mathcal N}$.
Let $\KK^{-1}_{af}({\mathcal N})$ denote the almost flat elements of $\KK^{-1}({\mathcal N})$.

\begin{definition} \label{2.5}
Let $\Gamma$ be a discrete group.
Given $\eta \in \KK^{-1}(B\Gamma)$, we say that $\eta$ is almost flat if for every 
compact manifold-with-boundary ${\mathcal N}$ and every
continuous map $\nu : {\mathcal N} \rightarrow B\Gamma$, the pullback $\nu^* \eta$ is almost flat.
\end{definition}

Let $\KK^{-1}_{af}(B\Gamma)$ denote the almost flat classes. It is a subgroup of $\KK^{-1}(B\Gamma)$.

The generator $\beta_{S^1}$ of
    $\KK^{-1}(S^1) \cong \Z$ is almost flat \cite[Exemple 3.2]{Skandalis (1991)}.
Let $V_{S^1}$, $\nabla^{V_{S_1}}$ and $\sigma^{V_{S_1}}$ be as in Definition \ref{2.4}, for 
$\beta_{S^1}$.

We will describe a map $\zeta : \KK^{*}_{af}({\mathcal N}) \rightarrow
\KK^{*-1}_{af}({\mathcal N} \times S^1)$ with the property that
$\int_{S^1} \circ \ch \circ \zeta = \ch$ when acting on $\KK^{*}_{af}({\mathcal N})$.
Let $\pi_1 : {\mathcal N} \times S^1 \rightarrow {\mathcal N}$ and
$\pi_2 : {\mathcal N} \times S^1 \rightarrow S^1$ be the projection maps. 

Consider first the case
when $\star = 0$. Let $V_{\mathcal N}$ and $\nabla^\pm_{\mathcal N}$ be as in Definition \ref{2.1}. 
The class of $\zeta([V_{\mathcal N}])$ in $\KK^{*-1}({\mathcal N} \times S^1)$
is represented by the ungraded vector bundle $V_{{\mathcal N} \times S^1} = (\pi_1^*
V_{\mathcal N}^+ \otimes \pi_2^* V_{S^1}) \oplus (\pi_1^*
V_{\mathcal N}^- \otimes \pi_2^* V_{S^1})$ with the automorphism $(\Id_{\pi_1^* V_{\mathcal N}^+}  \otimes \pi_s^* \sigma^{V_{S_1}})
\oplus (\Id_{\pi_1^* V_{\mathcal N}^-}  \otimes \pi_2^* (\sigma^{V_{S_1}})^{-1})$. Considering the connection
$\pi_1^* (\nabla^+_{\mathcal N} \oplus \nabla^-_{\mathcal N}) \otimes \pi_2^* \nabla^{V_{S_1}}$
shows that the $\KK$-theory class is almost flat.

Next, suppose that $\star = 1$. Let $V_{\mathcal N}$, $\nabla_{\mathcal N}$ and $\sigma_{\mathcal N}$ be as in Definition \ref{2.4}. Then $\zeta([V_{\mathcal N}, \sigma_{\mathcal N}])$ is represented by
taking the pullback of $V_{\mathcal N}$ under the map ${\mathcal N} \times [0,1] \rightarrow {\mathcal N}$
and gluing the two ends with $\sigma_{\mathcal N}$, to obtain a vector bundle on ${\mathcal N} \times S^1$. 
The connection $dt \partial_t + t \nabla_{\mathcal N} + (1-t) \sigma_{\mathcal N} \circ \nabla_{\mathcal N} 
\circ (\sigma_{\mathcal N})^{-1}$ on the pullback vector bundle gives a connection on the glued vector bundle, which
shows that the $\KK$-theory class is almost flat.

Similarly, we obtain a map $\zeta : \KK^{*}_{af}(B\Gamma) \rightarrow
\KK^{*-1}_{af}(B\Gamma \times S^1)$.
    
    \begin{corollary} \label{2.6}
    Given
    \begin{itemize}
\item    A closed
spin manifold $M$ with a 
positive scalar curvature metric,
\item A discrete group $\Gamma$,
\item A continuous map $\nu : M \rightarrow B\Gamma$, and 
\item An element $\eta \in
\KK^*_{af}(B \Gamma)$,
\end{itemize}
we have $\int_M \widehat{A}(TM) \wedge
\ch(\nu^* \eta) = 0$.
    \end{corollary}

\begin{corollary} \label{2.7}
 Given
   a closed connected
spin manifold $M$ with a 
positive scalar curvature metric, suppose that 
$\KK^*_{af}(B \pi_1(M, m_0)) \otimes \Q = \KK^*(B \pi_1(M, m_0)) \otimes \Q$.
Let $[M] \in \HH_{\dim(M)}(M; \Q)$ be the fundamental class.
Let $\nu : M \rightarrow B\pi_1(M, m_0)$ be the classifying map.
Then $\nu_* [M]$ vanishes in $\HH_{\dim(M)}(B \pi_1(M, m_0)); \Q)$.
\end{corollary}
\begin{proof}
Given $\eta \in \KK^*(B \pi_1(M, m_0))$, we have
$\int_M \widehat{A}(TM) \wedge \ch(\nu^* \eta) = \langle \nu_*(\star \widehat{A}(TM)), \ch(\eta) \rangle$,
where $\star \widehat{A}(TM) \in \HH_*(M; \Q)$ is the Poincar\'e dual of $\widehat{A}(TM) \in \HH^*(M; \Q)$.
Under the assumptions of the corollary, $\nu_*(\star \widehat{A}(TM))$ vanishes in $\HH_*(B \pi_1(M, m_0); \Q)$.
As $\star(1) = [M]$, the corollary follows.
\end{proof}

\subsection{Boundary obstruction} \label{subsect2.2}

    The next result gives an obstruction for a manifold-with-boundary to have a metric with positive scalar curvature
    in a neighborhood of the boundary, with nonnegative mean curvature on the boundary.
    We will use the following running assumption in this subsection.

    \begin{assumption} \label{2.8}
    Given $r_0, D > 0$, put $r_0^\prime = \frac{1}{256} r_0^2 D^2$
and $D^\prime = D + \frac{32}{r_0 D}$. 
Let $M$ be a connected Riemannian spin manifold-with-boundary, possibly incomplete.  Suppose that $\partial M$ is compact and has nonnegative mean curvature.  Let $K$ be a compact subset of $M$ containing $\partial M$.  

Suppose that 
\begin{itemize}
\item The distance neighborhood $N_{D^\prime}(K)$ lies in a compact submanifold-with-boundary ${\mathcal C}$,
\item $R > 0$ on $K$,
\item $R \ge r_0$ on $N_{D}(K) - K$ and
\item $R \ge - r_0^\prime$ on $N_{D^\prime}(K) - N_{D}(K)$.
\end{itemize}
    \end{assumption}
    
\begin{theorem} \label{2.9}
Suppose that Assumption \ref{2.8} holds.  Given $\beta \in 
\KK^{-1}_{af}({\mathcal C})$ if $M$ is even dimensional, or
$\beta \in 
\KK^{0}_{af}({\mathcal C})$ if $M$ is odd dimensional, we have
\begin{equation} \label{2.10}
\int_{\partial M} \widehat{A}(T\partial M) \wedge \ch \left(
\beta \big|_{\partial M} \right) = 0.
\end{equation}
\end{theorem}
\begin{proof}
Suppose first that $M$ is even
dimensional.
Let $\epsilon > 0$ be a small parameter, which we will adjust.
Let $d_K \in C({\mathcal C})$ be the distance function from $K$. 
We would like to use level sets of $d_K$ but they need not be smooth.  To get around this, let
$\widetilde{d}_K$ be a slight smoothing of $d_K$ on $N_{D^\prime}(K)$ so that
$|\widetilde{d}_K - d_K| < \epsilon$ and $|\nabla \widetilde{d}_K | < 1 + \epsilon$. 
Put $\lambda_1 = D-\epsilon$. 
Choose a regular value $\lambda_2$ of $\widetilde{d}_K$ that
is slightly greater than $D-\epsilon + \frac{16}{r_0(D-2\epsilon)}$. 
Put $\lambda_0 = \lambda_1 - 
\frac{16}{r_0(\lambda_2 - \lambda_1)}$. Then $\lambda_0$ is slightly greater than $\epsilon$. By construction,
\begin{align} \label{2.11}
\widetilde{d}_K^{-1}(- \infty, \lambda_0)
& \subset N_{3\epsilon}(K), \\
\widetilde{d}_K^{-1}(\lambda_0, \lambda_1)
& \subset N_D(K) - K, \notag \\
\widetilde{d}_K^{-1}(\lambda_1, \lambda_2)
& \subset N_{D^\prime}(K). \notag
\end{align}
In particular,
\begin{equation} \label{2.12} 
\begin {cases}
R(x) > 0 & \text{if $\widetilde{d}_K(x) \le \lambda_0$,} \\
R(x) \ge r_0 & \text{if $\lambda_0 \le \widetilde{d}_K(x) \le \lambda_1$,} \\
R(x) \ge - r_0^\prime & \text{if $\lambda_1 \le  \widetilde{d}_K(x)< \lambda_2$.}
\end{cases}
\end{equation}

Define $\sigma \in C(- \infty, \lambda_2)$ by
\begin{equation} \label{2.13}
\sigma(t) = 
\begin{cases}
0 & \text{if $t \le \lambda_0$,} \\
\frac{r_0}{8} (t - \lambda_0) & \text{if $\lambda_0 \le t \le \lambda_1$,} \\
\frac{2}{\lambda_2 - t} & \text{if $\lambda_1 \le t < \lambda_2$.}
\end{cases}
\end{equation}
Let $\widetilde{\sigma}$ be a slight smoothing of $\sigma$.
Put $F = \sigma \circ \widetilde{d}_K$ and 
$f = \widetilde{\sigma} \circ \widetilde{d}_K$.

The preimage $\widetilde{d}_K^{-1}(- \infty, \lambda_2]$ is a smooth
compact manifold-with-boundary. Its boundary is the disjoint union of 
$\partial M$ and $\widetilde{d}_K^{-1}(\lambda_2)$. For small $\epsilon$,
let
${\mathcal N}^\epsilon$ be the codimension-zero submanifold of $\widetilde{d}_K^{-1}(- \infty, \lambda_2]$ whose
boundary consists of $\partial M$ along with interior points of distance
$\epsilon$ from $\widetilde{d}_K^{-1}(\lambda_2)$. We write
$\partial {\mathcal N}^\epsilon = \partial M \cup \partial_+{\mathcal N}^\epsilon$.

Let $\beta^\epsilon \in 
\KK^{-1}_{af}({\mathcal N}^\epsilon)$ be the restriction of $\beta$ to 
${\mathcal N}^\epsilon$. As in
Definition \ref{2.4}, 
let $(V, \sigma)$ be a Hermitian vector bundle on ${\mathcal N}^\epsilon$ 
with unitary automorphism that represents $\beta^\epsilon$.
Let $E$ be the $\Z_2$-graded vector bundle $V \oplus V$. The operator
$\begin{pmatrix}
0 & \sigma^{-1} \\
\sigma & 0
\end{pmatrix}$ is an odd self-adjoint endomorphism of $E$. Put 
\begin{equation} \label{2.14}
A = 
\begin{pmatrix}
\nabla^V & f \sigma^{-1} \\
f\sigma & \nabla^V
\end{pmatrix},
\end{equation}
a superconnection on $E$. 
Let $D^E$ be the quantization of $A$, i.e.
\begin{equation} \label{2.15}
D^E = 
\begin{pmatrix}
D^V & \epsilon_S f \sigma^{-1} \\
\epsilon_S f \sigma & D^V
\end{pmatrix},
\end{equation}
where $\epsilon_S$ is the $\Z_2$-grading
operator on the spinor bundle $S$ and $D^V$ is the Dirac operator on $C^\infty(M; S \otimes V)$. Then
\begin{align} \label{2.16}
(D^E)^2 = 
& \begin{pmatrix}
(D^V)^2 + f^2 & 0 \\
0 & (D^V)^2 + f^2
\end{pmatrix} +  \\
& \sqrt{-1} \epsilon_S c(df)
\begin{pmatrix}
0 &  \sigma^{-1} \\
 \sigma & 0
\end{pmatrix} + \sqrt{-1} \epsilon_S f
\begin{pmatrix}
0 & - \sigma^{-1} c(\nabla^V \sigma) \sigma^{-1} \\
 c(\nabla^V \sigma) & 0
\end{pmatrix}, \notag
\end{align}
where $c(w)$ denotes Clifford multiplication by a $1$-form $w$ (satisfying $c(w)^2 = |w|^2$). 

For notation, if $\{e_\alpha\}_{\alpha = 1}^n$ is a
local orthonormal frame of $TM$ and $\{\tau^\alpha\}_{\alpha = 1}^n$
is the dual coframe then we put $\gamma^\alpha = 
c(\tau^\alpha)$.  Locally, $D^V = - \sqrt{-1}
\sum_{\alpha = 1}^n \gamma^\alpha \nabla_{e_\alpha}$.

We take the inner
product $\langle \psi_1, \psi_2 \rangle$ on $S \otimes E$ to be $\C$-linear in $\psi_2$ and $\C$-antilinear in $\psi_1$.
On $\partial {\mathcal N}^\epsilon$, we let $e_n$ be the inward pointing
normal vector. 
For $\psi_1, \psi _2 \in C^\infty({\mathcal N}^\epsilon;
S \otimes E)$, we have
\begin{align} \label{2.17}
& \int_{{\mathcal N}^\epsilon} \langle D^E \psi_1, \psi_2 \rangle  \; \dvol_{{\mathcal N}^\epsilon} -
\int_{{\mathcal N}^\epsilon} \langle \psi_1, D^E \psi_2 \rangle  \; \dvol_{{\mathcal N}^\epsilon} = \\
& - \sqrt{-1} \int_{\partial {\mathcal N}^\epsilon} \langle \psi_1, \gamma^n \psi_2 \rangle  \; \dvol_{\partial {\mathcal N}^\epsilon}. \notag
\end{align}

Define a self-adjoint operator $\Pi$ on $C^\infty \left( \partial 
{\mathcal N}^\epsilon, (S \otimes E) \big|_{\partial {\mathcal N}^\epsilon} \right)$ by
\begin{equation} \label{2.18}
\Pi = 
\begin{cases}
\begin{pmatrix}
0 & \sqrt{-1} \epsilon_S \gamma^n \\
\sqrt{-1} \epsilon_S \gamma^n & 0
\end{pmatrix} & \text{on $\partial M$}, \\
\begin{pmatrix}
0 & \sqrt{-1} \epsilon_S \gamma^n \sigma^{-1} \\
\sqrt{-1} \epsilon_S \gamma^n \sigma & 0
\end{pmatrix} & \text{on $\partial_+{\mathcal N}^\epsilon$}
\end{cases}.
\end{equation}
Then $\Pi^2 = \Id$, $\Pi \gamma^n + \gamma^n \Pi = 0$ and $\Pi$ commutes with the $\Z_2$ grading on
$(S \otimes E) \big|_{\partial {\mathcal N}^\epsilon}$. It follows
from (\ref{2.17}) that with the boundary condition $\Pi \psi = \psi$,
the operator $D^E$ on $C^\infty({\mathcal N}^\epsilon;
S \otimes E)$ is formally self-adjoint. In fact, the boundary condition is
an elliptic boundary condition \cite[Section 7.5]{Baer-Ballmann (2012)}.
We let ${\mathcal D}$ denote the ensuing self-adjoint operator, densely
defined on $L^2({\mathcal N}^\epsilon;
S \otimes E)$.

\begin{lemma} \label{2.19}
There is an $\epsilon_0 > 0$ so that if $\epsilon < \epsilon_0$, and in addition
$\| F^V \| \le \epsilon$ and $\| \nabla^V \sigma \| \le \epsilon$, then
the kernel of ${\mathcal D}$ vanishes.
    \end{lemma}
    \begin{proof}
    One can check that
\begin{equation} \label{2.20} 
|\nabla F|(x) =   
\begin {cases}
0 & \text{if $\widetilde{d}_K(x) < \lambda_0$,} \\
\frac{r_0}{8} |\nabla \widetilde{d}_K|(x) & \text{if $\lambda_0 < \widetilde{d}_K(x) < \lambda_1$,} \\
\frac{2}{(\lambda_2 - \widetilde{d}_K(x))^2} |\nabla \widetilde{d}_K|(x)  & \text{if $\lambda_1 <  \widetilde{d}_K(x) < \lambda_2$.}
\end{cases}
\end{equation}
Then
\begin{equation} \label{2.21} 
|\nabla F|(x) \le   
\begin {cases}
0 & \text{if $\widetilde{d}_K(x) < \lambda_0$,} \\
\frac{r_0}{8} (1+\epsilon) & \text{if $\lambda_0 < \widetilde{d}_K(x) < \lambda_1$,} \\
\frac{2}{(\lambda_2 - \widetilde{d}_K(x))^2} (1+\epsilon)  & \text{if $\lambda_1 <  \widetilde{d}_K(x) < \lambda_2$.}
\end{cases}
\end{equation}
Writing $F^2 - \epsilon F  = \left( F -\frac{\epsilon}{2} \right)^2  - \frac{\epsilon^2}{4}$, we have
\begin{equation} \label{2.22}
 (F^2 - \epsilon F)(x) \ge   
   \begin {cases}
0 & \text{if $\widetilde{d}_K(x) < \lambda_0$,} \\
 - \frac{\epsilon^2}{4} & \text{if $\lambda_0 < \widetilde{d}_K(x) < \lambda_1$,} \\
\left( \frac{2}{\lambda_2 - \widetilde{d}_K(x)} -\frac{\epsilon}{2} \right)^2 - \frac{\epsilon^2}{4} & \text{if $\lambda_1 <  \widetilde{d}_K(x) < \lambda_2$.}
\end{cases}     
\end{equation}
For small $\epsilon$, we can bound $\left( \frac{2}{\lambda_2 - \widetilde{d}_K(x)} -\frac{\epsilon}{2} \right)^2 - \frac{2}{(\lambda_2 - \widetilde{d}_K(x))^2} (1+\epsilon)$ below by $\frac{1}{(\lambda_2 - \widetilde{d}_K(x))^2}$. Using the fact that $\lambda_2 - \lambda_1 = \frac{16}{r_0 (D-2\epsilon)}$, it follows that for
small $\epsilon$, we have
\begin{equation} \label{2.23}
 (F^2 - |\nabla F| - \epsilon F)(x) \ge   
   \begin {cases}
0 & \text{if $\widetilde{d}_K(x) < \lambda_0$,} \\
 - \frac{\epsilon^2}{4} - \frac{r_0}{8} (1+\epsilon) & \text{if $\lambda_0 < \widetilde{d}_K(x) < \lambda_1$,} \\
\frac{r_0^2 D^2}{512} & \text{if $\lambda_1 <  \widetilde{d}_K(x) < \lambda_2$.}
\end{cases}    
\end{equation}
Combining with (\ref{2.12}), we see that for small $\epsilon$, the function
$\frac{R}{4} + F^2 - |\nabla F| - \epsilon F$ is bounded below on 
its domain in $\widetilde{d}_K^{-1}(- \infty, \lambda_2]$ by a positive constant that
is independent of $\epsilon$.
Taking $\widetilde{\sigma}$ to be an appropriate smoothing of
$\sigma$, the same will be true of
$\frac{R}{4} + f^2 - |\nabla f| - \epsilon f$.

The Lichnerowicz formula gives
\begin{equation} \label{2.24}
(D^V)^2 = \nabla^* \nabla + \frac{R}{4} - \frac{1}{4} [\gamma^\mu, \gamma^\nu]
F^V_{\mu \nu}.
\end{equation}
If $\psi \in C^\infty({\mathcal N}^\epsilon;
S \otimes E)$ then
\begin{equation} \label{2.25}
\int_{{\mathcal N}^\epsilon} \langle \psi, \nabla^* \nabla \psi
\rangle \: \dvol_{{\mathcal N}^\epsilon} = 
\int_{{\mathcal N}^\epsilon} | \nabla \psi |^2
 \: \dvol_{{\mathcal N}^\epsilon} + 
\int_{\partial {\mathcal N}^\epsilon} \langle \psi,
\nabla_{e_n} \psi \rangle \: \dvol_{\partial {\mathcal N}^\epsilon}.
\end{equation}
Suppose that $\psi$ is a nonzero solution of $D^E \psi = 0$. Then
$\int_{{\mathcal N}^\epsilon} \langle \psi, (D^E)^2 \psi \rangle \: 
\dvol_{{\mathcal N}^\epsilon} = 0$, which from (\ref{2.16}) becomes
\begin{align} \label{2.26}
 0 = & 
\int_{{\mathcal N}^\epsilon} \left( | \nabla \psi |^2 + 
\left( \frac{R}{4} + f^2 \right) |\psi|^2 - \frac14 
\langle \psi, [\gamma^\alpha, \gamma^\beta] F^V_{\alpha \beta} \psi
\rangle \right) \: \dvol_{{\mathcal N}^\epsilon}
+ \\
&\int_{{\mathcal N}^\epsilon}
\left\langle \psi, 
\sqrt{-1} \epsilon_S c(df)
\begin{pmatrix}
0 &  \sigma^{-1} \\
 \sigma & 0
\end{pmatrix} \psi \right\rangle \dvol_{{\mathcal N}^\epsilon}
+ \notag \\
& \int_{{\mathcal N}^\epsilon}
\left\langle \psi,
\sqrt{-1} \epsilon_S f
\begin{pmatrix}
0 & - \sigma^{-1} c(\nabla^V \sigma) \sigma^{-1} \\
 c(\nabla^V \sigma) & 0
\end{pmatrix} \psi \right\rangle \dvol_{{\mathcal N}^\epsilon}  + \notag \\
& \int_{\partial {\mathcal N}^\epsilon} \langle \psi,
\nabla_{e_n} \psi \rangle \: \dvol_{\partial {\mathcal N}^\epsilon}. \notag
\end{align}
From $D^E \psi = 0$, one can solve for $\nabla_{e_n} \psi$ and obtain
\begin{align} \label{2.27}
\int_{\partial {\mathcal N}^\epsilon} \langle \psi,
\nabla_{e_n} \psi \rangle \: \dvol_{\partial {\mathcal N}^\epsilon} =
& \int_{\partial {\mathcal N}^\epsilon} \langle \psi,
D^E_{\partial {\mathcal N}^\epsilon} \psi \rangle \: \dvol_{\partial {\mathcal N}^\epsilon} + 
\int_{\partial {\mathcal N}^\epsilon} H | \psi |^2 \: \dvol_{\partial {\mathcal N}^\epsilon} + \\
& \int_{\partial {\mathcal N}^\epsilon} \left\langle \psi,
\begin{pmatrix}
0 & \sqrt{-1} f \epsilon_S \gamma^n \sigma^{-1} \\
\sqrt{-1} f \epsilon_S \gamma^n \sigma & 0
\end{pmatrix}
\psi \right\rangle \: \dvol_{\partial {\mathcal N}^\epsilon}, \notag
\end{align}
where 
\begin{equation} \label{2.28}
D^E_{\partial {\mathcal N}^\epsilon} = - \sum_{j=1}^{n-1} \gamma^n
\gamma^j \nabla_{e_j} \otimes I_2
\end{equation}
is the intrinsic Dirac
operator on ${\partial {\mathcal N}^\epsilon}$ coupled to
$E \Big|_{\partial {\mathcal N}^\epsilon}$.

We now impose the boundary condition $\Pi \psi = \psi$. 
Then (\ref{2.27}) becomes
\begin{equation} \label{2.29}
\int_{\partial {\mathcal N}^\epsilon} \langle \psi,
\nabla_{e_n} \psi \rangle \: \dvol_{\partial {\mathcal N}^\epsilon} =
\int_{\partial {\mathcal N}^\epsilon} \langle \psi,
D^E_{\partial {\mathcal N}^\epsilon} \psi \rangle \: \dvol_{\partial {\mathcal N}^\epsilon} + 
\int_{\partial {\mathcal N}^\epsilon} (H+f) | \psi |^2 \: \dvol_{\partial {\mathcal N}^\epsilon}.
\end{equation}
On
$\partial M$, we have $D^E_{\partial {\mathcal N}^\epsilon} \Pi + 
\Pi D^E_{\partial {\mathcal N}^\epsilon} = 0$, which implies that
$\int_{\partial M} \langle \psi,
D^E_{\partial {\mathcal N}^\epsilon} \psi \rangle \: \dvol_{\partial M} = 0$. On $\partial_+ {\mathcal N}^\epsilon$ we have
\begin{equation} \label{2.30}
D^E_{\partial {\mathcal N}^\epsilon} \Pi + 
\Pi D^E_{\partial {\mathcal N}^\epsilon} =
\begin{pmatrix}
0 & - \sqrt{-1} \epsilon_S \gamma^j \sigma^{-1} (\nabla^V_{e_j} \sigma) 
\sigma^{-1} \\
\sqrt{-1} \epsilon_S \gamma^j \nabla^V_{e_j} \sigma & 0
\end{pmatrix},
\end{equation}
which gives
\begin{align} \label{2.31}
& \int_{\partial_+ {\mathcal N}^\epsilon} \langle \psi,
D^E_{\partial {\mathcal N}^\epsilon} \psi \rangle \: 
\dvol_{\partial_+ {\mathcal N}^\epsilon} = 
\frac12 \int_{\partial_+ {\mathcal N}^\epsilon} \langle \psi,
\left( D^E_{\partial {\mathcal N}^\epsilon} \Pi + \Pi D^E_{\partial {\mathcal N}^\epsilon} \right) \psi \rangle \: 
\dvol_{\partial_+ {\mathcal N}^\epsilon} =
\\
& \frac12 \int_{\partial_+ {\mathcal N}^\epsilon} \left\langle \psi,
\begin{pmatrix}
0 & - \sqrt{-1} \epsilon_S \gamma^j \sigma^{-1} (\nabla^V_{e_j} \sigma) 
\sigma^{-1} \\
\sqrt{-1} \epsilon_S \gamma^j \nabla^V_{e_j} \sigma & 0
\end{pmatrix}
\psi \right\rangle \: 
\dvol_{\partial_+ {\mathcal N}^\epsilon}. \notag
\end{align}
Hence
\begin{align} \label{2.32}
& \int_{\partial {\mathcal N}^\epsilon} \langle \psi,
\nabla_{e_n} \psi \rangle \: 
\dvol_{\partial {\mathcal N}^\epsilon} \ge \\
& \int_{\partial M} H |\psi|^2 \dvol_{\partial M} +
\int_{\partial_+ {\mathcal N}^\epsilon}
(H + f - \const \epsilon) |\psi|^2 \: \dvol_{\partial_+ {\mathcal N}^\epsilon}, \notag
\end{align}
where $\const = \const(n)$.

Referring to the terms in (\ref{2.26}), we have a bound
\begin{align} \label{2.33}
& \int_{{\mathcal N}^\epsilon} 
\left( \frac{R}{4} + f^2 \right) |\psi|^2 
 \: \dvol_{{\mathcal N}^\epsilon}
+ \int_{{\mathcal N}^\epsilon}
\left\langle \psi, 
\sqrt{-1} \epsilon_S c(df)
\begin{pmatrix}
0 &  \sigma^{-1} \\
 \sigma & 0
\end{pmatrix} \psi \right\rangle \dvol_{{\mathcal N}^\epsilon} + \\
& \int_{{\mathcal N}^\epsilon}
\left\langle \psi,
\sqrt{-1} \epsilon_S f
\begin{pmatrix}
0 & - \sigma^{-1} c(\nabla^V \sigma) \sigma^{-1} \\
 c(\nabla^V \sigma) & 0
\end{pmatrix} \psi \right\rangle \dvol_{{\mathcal N}^\epsilon}  \ge \notag \\
& \int_{{\mathcal N}^\epsilon} 
\left( \frac{R}{4} + f^2 - |\nabla f| - \epsilon f \right) |\psi|^2 
 \: \dvol_{{\mathcal N}^\epsilon}. \notag
\end{align}
Since $\widetilde{d}_K^{-1}(\lambda_2)$ is a smooth compact hypersurface in
$M$, for small $\epsilon$ there is a uniform bound on the mean curvature
$H$ of $\partial_+ {\mathcal N}^\epsilon$, independent of $\epsilon$.
As $f$ blows up near $\widetilde{d}_K^{-1}(\lambda_2)$, equation
(\ref{2.32}) implies that for small
$\epsilon$ the expression 
$\int_{\partial {\mathcal N}^\epsilon} \langle \psi,
\nabla_{e_n} \psi \rangle \: \dvol_{\partial {\mathcal N}^\epsilon}$ is
nonnegative. Finally, 
\begin{equation} \label{2.34}
-  \frac14 
\int_{{\mathcal N}^\epsilon}  
\langle \psi, [\gamma^\alpha, \gamma^\beta] F^V_{\alpha \beta} \psi
\rangle \: \dvol_{{\mathcal N}^\epsilon} \ge - \const
\epsilon \int_{{\mathcal N}^\epsilon}  
|\psi|^2 \: \dvol_{{\mathcal N}^\epsilon}.
\end{equation}
Using the uniform positivity of $\frac{R}{4} + f^2 - |\nabla f| - \epsilon f$ for small $\epsilon$, we see that if $\epsilon$ is sufficiently small
then we obtain a contradiction to (\ref{2.26}). This proves the lemma.
    \end{proof}

We will loosely refer to the index of ${\mathcal D}$, meaning the index of the restriction of
${\mathcal D}$ to $\Dom({\mathcal D})^+$.

    \begin{lemma} \label{2.35}
    The index of ${\mathcal D}$ equals $\int_{\partial M} \widehat{A}(T\partial M) \wedge 
    \ch \left( \beta \big|_{\partial M} \right)$.
    \end{lemma}
\begin{proof}
As we have local boundary conditions, the index only depends on the principal
symbols of $D^E$ and $\Pi$, up to homotopy. We can perform deformations of the geometry without changing the index, as long as we preserve ellipticity.

We first make the Riemannian metric on ${\mathcal N}^\epsilon$ a
product near $\partial {\mathcal N}^\epsilon$. 
Let us then parametrize a product neighborhood of $\partial {\mathcal N}^\epsilon$ as
$\partial {\mathcal N}^\epsilon \times [0, \delta)$ for some small $\delta > 0$.
Next, we deform the
function $f$ to zero.  
Deform $\sigma$ on $\partial {\mathcal N}^\epsilon \times [0, 3\delta/4)$ so that
it pulls back from $\partial {\mathcal N}^\epsilon$ there.
Deform the connection $\nabla^V$ on $V$
to be a product on $\partial {\mathcal N}^\epsilon \times [0, 3\delta/4)$.
Construct a connection $\nabla^E$ on $E = V \oplus V$ so that
on $\partial M \times [0, \delta/4]$, we have
$\nabla^E = 
\begin{pmatrix}
\nabla^V & 0 \\
0 & \nabla^V
\end{pmatrix}$, and on the complement of $\partial M \times [0, \delta/2]$, 
we have
$\nabla^E = 
\begin{pmatrix}
\nabla^V & 0 \\
0 & \sigma \circ \nabla^V \circ \sigma^{-1}
\end{pmatrix}$. Then 
$\begin{pmatrix}
0 & \sigma^{-1} \\
\sigma & 0    
\end{pmatrix} 
\nabla^E
\begin{pmatrix}
0 & \sigma^{-1} \\
\sigma & 0    
\end{pmatrix}
= \nabla^E$ on the complement of $\partial M \times [0, \delta/2]$.

We now use an argument as in \cite[Proof of Theorem 1.1]{Lott (2021)}.
Let $D{\mathcal N}^\epsilon$ denote the double of ${\mathcal N}^\epsilon$,
which we can consider to be the result of gluing two disjoint copies ${\mathcal N}_1^\epsilon$ and
${\mathcal N}_2^\epsilon$ of ${\mathcal N}^\epsilon$, along their boundaries.
Let $E_1$ and $E_2$ be copies of $E$, on ${\mathcal N}_1^\epsilon$ and
${\mathcal N}_2^\epsilon$, respectively. Let
$\nabla^{E_1}$ and $\nabla^{E_2}$ be copies of $\nabla^E$ acting on 
sections of $E_1$ and $E_2$, respectively.
Let $DE$ be the vector bundle on $D{\mathcal N}^\epsilon$ obtained by gluing
$E_1$ and $E_2$ along $\partial M$ by $\begin{pmatrix}
0 & 1 \\
1 & 0    
\end{pmatrix}$, and along $\partial_+ {\mathcal N}^\epsilon$ by
$\begin{pmatrix}
0 & \sigma^{-1} \\
\sigma & 0    
\end{pmatrix}$.
Then $\nabla^{E_1}$ and $\nabla^{E_2}$ descend to give a connection $\nabla^{DE}$ on $DE$. 

The involution on the disjoint union $E_1 \coprod E_2$, that sends an element of $E_1$ to the 
corresponding element of $E_2$ and vice versa, descends to an involution of $DE$ that covers
the involution on $D{\mathcal N}^\epsilon$. The involution acts on
$DE \big|_{\partial M}$ as $\begin{pmatrix}
0 & 1 \\
1 & 0    
\end{pmatrix}$, and on $DE \big|_{\partial_+ {\mathcal N}^\epsilon}$ as
$\begin{pmatrix}
0 & \sigma^{-1} \\
\sigma & 0    
\end{pmatrix}$.
It preserves the connection $\nabla^{DE}$.

Then the boundary condition $\Pi \psi = \psi$ amounts to
looking at $\Z_2$-invariant sections of $DS \otimes DE$. Hence
the index of ${\mathcal D}$ is the same as the index of the
Dirac-type operator acting on $\Z_2$-invariant sections of $DS \otimes DE$,
which is same as the index of the orbifold Dirac-type operator on
the $\Z_2$-quotient of $D{\mathcal N}^\epsilon$. Since
$\partial {\mathcal N}^\epsilon$ is odd dimensional, it does not
contribute to the orbifold index formula.  Hence
the index of ${\mathcal D}$ equals $\int_{{\mathcal N}^\epsilon}
\widehat{A}(T{\mathcal N}^\epsilon) \wedge \tr_s \left( e^{ \frac{i}{2\pi} (\nabla^E)^2}
\right)$. By construction,
\begin{align} \label{2.36} 
\int_{{\mathcal N}^\epsilon}
\widehat{A}(T{\mathcal N}^\epsilon) \wedge \tr_s \left( e^{  \frac{i}{2\pi}(\nabla^E)^2}
\right) = &
\int_{\partial M \times [0, \delta]}
\widehat{A}(T{\mathcal N}^\epsilon) \wedge 
\tr_s \left( e^{ \frac{i}{2\pi} (\nabla^E)^2}
\right) \\ 
= &
\int_{\partial M}
\widehat{A}(T \partial M) \wedge 
\int_{[0, \delta]}
\tr_s \left( e^{ \frac{i}{2\pi} (\nabla^E)^2} \right) \notag \\
= & 
\int_{\partial M}
\widehat{A}(T\partial M) \wedge \ch(V, \sigma) \big|_{\partial M}  \notag \\
= & 
\int_{\partial M} \widehat{A}(T\partial M) \wedge \ch \left( \beta \big|_{\partial M}
\right).  \notag
\end{align}
This proves the lemma.
\end{proof}

    Combining Lemmas \ref{2.19} and \ref{2.35} proves the theorem when $M$ is even dimensional. Suppose that $M$ is odd dimensional.
    Consider $M^\prime = M \times S^1$ and $K^\prime = K \times S^1$.
    The theorem for $M$ now follows from the theorem for $M^\prime$.
\end{proof}

\begin{corollary} \label{2.38}
Suppose that Assumption \ref{2.8} holds.
Let $\Gamma$ be a discrete group and let $\eta$ be an element of $\KK^*_{af}(B\Gamma)$. Given a continuous map $\nu :  {\mathcal C} \rightarrow B\Gamma$, let $\nu \big|_{\partial M}$ be its restriction to $\partial M$.  Then
\begin{equation} \label{2.39}
\int_{\partial M} \widehat{A}(T\partial M) \wedge \ch \left( \nu \big|_{\partial M}^* \eta
\right) = 0.
\end{equation}
\end{corollary}

\begin{remark} \label{2.40}
There is some overlap between Corollary \ref{2.38} and
\cite[Theorem 6.12]{Gromov-Lawson (1983)}. Namely, if $M$ is a complete noncompact spin manifold-with-boundary with compact
boundary
then by taking ${\mathcal C}$ large enough, Corollary \ref{2.38} gives obstructions to $M$ having a complete Riemannian metric
with positive scalar curvature and nonnegative mean curvature boundary.  On the other hand, if there is such a metric then
taking the double $DM$ and smoothing the metric, we get a complete Riemannian metric with positive scalar curvature, to which
\cite[Theorem 6.12]{Gromov-Lawson (1983)} gives some obstructions.
\end{remark}

\subsection{End obstructions} \label{subsect2.3}

We show how Theorem \ref{2.9} gives end obstructions to the existence of complete finite volume metrics with
positive scalar curvature.
We first prove a general result about the existence of hypersurfaces with nonnegative mean curvature, as seen from 
infinity, in a finite volume manifold. The proof was explained to me by Antoine Song.

\begin{proposition} \label{2.41} 
Let $M$ be a complete connected
finite volume oriented Riemannian manifold with compact boundary, of dimension at most seven. Then there is an exhaustion of $M$ by compact connected submanifolds-with-boundary $Z$ so that $\partial Z$ has nonnegative mean curvature, as
seen from $M-Z$.
\end{proposition}
\begin{proof}
The proof is along the lines of that of \cite[Claim 2.4]{Song (2023)}.
Fixing a basepoint $m_0 \in M$, slightly smoothing the distance function from $m_0$ on compact subsets, and applying the coarea formula, we see that there is an exhaustion of $M$ by compact codimension-zero submanifolds $\{K_i\}_{i=1}^\infty$
with smooth boundary so that $\lim_{i \rightarrow \infty} \vol_{n-1}(\partial K_i) = 0$.  
The idea is now to move each $\partial K_i$ in a $1$-parameter family of hypersurfaces
within $K_i$. 
Given $i$, consider a $1$-parameter family
$L = \{L_t\}_{t \in [0,1]}$ of $n$-dimensional integral currents in $K_i$
so that 
\begin{itemize}
\item $\partial L_0 = \partial K_i$,
\item $\partial L_t$ is a continuous function of $t$ in the flat topology, and
\item For all $t$, $\vol_{n-1}(\partial L_t) \le 2 
\vol_{n-1}(\partial K_i)$.
\end{itemize}
Let $\Gamma_i$ denote the collection of such $1$-parameter families.

Fix a small ball $B$ around $m_0$.
We claim that for all large $i$, $\partial L_1$ is nonempty for any
$L \in \Gamma_i$.  Otherwise, the family 
$\{\partial L_t\}_{t \in [0,1]}$ would sweep out $K_i$.  In particular,
some $\partial L_t$ would divide $B$ into two parts of equal volume.
However, for large $i$,
this contradicts the assumption that $\vol_{n-1}(\partial L_t) \le 2 
\vol_{n-1}(\partial K_i)$.

Now take a sequence $\{L^j\}_{j=1}^\infty$ in $\Gamma_i$
so that $\vol_{n-1}(\partial L^j_1)$ approaches
$\inf_{L \in \Gamma_i} \vol_{n-1}(\partial L_1)$. After passing to a
subsequence, we can assume that $\lim_{j \rightarrow \infty} \partial L^j
= W_i$ in the flat topology, for some $C^{1,1}$-regular hypersurface
$W_i$.  By construction, $\vol_{n-1} W_i \le 2 
\vol_{n-1}(\partial K_i)$. Where $W_i$ doesn't intersect $\partial K_i$, it will be a minimal hypersurface.  Some of $W_i$ may be hung up on $\partial K_i$ but in any case, 
$W_i$ will have nonnegative mean curvature in the weak sense
as seen from $M-K_i$. Let $\widehat{Z}_i$ be the closure of the connected component of $M - W_i$
that contains $B$.  By an argument as in the proof of the above claim, the
$\widehat{Z}_i$'s exhaust $M$. Finally, by running the mean curvature
flow for a short time, we can smooth $W_i$ to a nearby smooth hypersurface
$H_i$ with nonnegative mean curvature.
Letting ${Z}_i$ denote the corresponding submanifold of $M$ with
boundary component $H_i$, we obtain the desired exhaustion of $M$.
\end{proof}

\begin{corollary} \label{2.42}
Let $M$ be a connected complete finite volume Riemannian spin manifold-with-boundary, with compact boundary, positive scalar curvature outside of a compact set, and dimension
at most seven.

Let $\Gamma$ be a discrete group and let $\eta$ be an element of $\KK^*_{af}(B\Gamma)$. Given a continuous map $\nu :  M \rightarrow B\Gamma$, let $\nu \big|_{\partial M}$ be its restriction to $\partial M$.  Then
\begin{equation} \label{2.43}
\int_{\partial M} \widehat{A}(T\partial M) \wedge \ch \left( \nu \big|_{\partial M}^* \eta
\right) = 0.
\end{equation}
\end{corollary}
\begin{proof}
We apply Proposition \ref{2.41} to get $Z$, with $\partial M \subset Z$ and
$\overline{M - Z}$ having positive scalar curvature. Let us write $\partial_+Z$ for $\partial \overline{M - Z}$, so
$\partial Z = \partial M \cup \partial_+ Z$. Applying
Theorem \ref{2.9} to a large compact submanifold ${\mathcal C}$ of $\overline{M-Z}$ containing $\partial_+Z$ gives
\begin{equation} \label{2.44}
\int_{\partial_+ Z} \widehat{A}(T\partial_+ Z) \wedge \ch \left( \nu \big|_{\partial_+ Z}^* \eta
\right) = 0.
\end{equation}
Deforming the metric on $Z$ to be a product near the boundary gives
\begin{align} \label{2.45}
& \int_{\partial_+ Z} \widehat{A}(T\partial_+ Z) \wedge \ch \left( \nu \big|_{\partial_+ Z}^* \eta
\right) - \int_{\partial M} \widehat{A}(T\partial M) \wedge \ch \left( \nu \big|_{\partial M}^* \eta
\right) = \\
& \int_Z d
\left( \widehat{A}(TZ) \wedge \ch \left( \nu \big|_{Z}^* \eta
\right) \right) = 0. \notag
\end{align}
Equation (\ref{2.43}) follows.
\end{proof}

\begin{example} \label{2.46}
Let $C$ be a compact connected spin manifold of dimension at most six.  Let $\nu : C \rightarrow B\pi_1(C)$ be
the classifying map for the universal cover of $C$, defined up to
homotopy.  If $[0, \infty) \times C$ has a complete finite volume
Riemannian metric with positive scalar curvature then for any $\eta \in \KK^*_{af}(B \pi_1(C))$, we have
$\int_C \widehat{A}(TC) \wedge \ch \left( \nu^* \eta \right) = 0$.

For example, there is no finite volume complete Riemannian metric with positive scalar curvature on $[0, \infty) \times T^{n-1}$ if
$n \le 7$.  Of course, such metrics exist if
one removes the finite volume condition.
\end{example}

\begin{remark} \label{2.47}
To see more clearly that Theorem \ref{2.9} gives an end obstruction to
positive scalar curvature, suppose that $M$ is a complete noncompact
finite volume Riemannian spin manifold having $\dim(M) \le 7$, with positive
scalar curvature outside of a compact subset. For notational simplicity,
suppose that $\dim(M)$ is even; the case when $\dim(M)$ is odd is similar.
Let $\{K_i\}$ be an
exhaustion of $M$ by compact submanifolds-with-boundary, with
$K_i \subset \Int(K_{i+1})$.
Put ${\mathcal K}^{-1}_{af,i} = \lim_{\stackrel{j}{\leftarrow}} \KK^{-1}_{af}(K_j - 
\Int(K_i))$, an
inverse limit, and
then put ${\mathcal K}^{-1}_{af, \infty} = \lim_{\stackrel{i}{\rightarrow}} {\mathcal K}^{-1}_{af,i}$,
a direct limit. Clearly ${\mathcal K}^{-1}_{af, \infty}$ only depends on the geometry of $M$ at
infinity. 
Given $\beta \in {\mathcal K}^{-1}_{af, \infty}$, we obtain an obstruction
as follows.  Choose a $\beta_i \in {\mathcal K}^{-1}_{af,i}$ representing $\beta$, for some $i$.
Taking $i$ large, we can assume that $M-K_i$ has positive scalar curvature.
Take $\beta_{i,j} \in 
\KK^{-1}_{af}(K_j -  \Int(K_i))$ in the inverse limit defining $\beta_i$. 
Spin cobordism implies that $\int_{\partial K_i} \widehat{A}(T\partial K_i) \wedge \ch \left(\beta_{i,j} \big|_{\partial K_i}
\right)$ is independent of the choices made. Taking $j$ sufficiently large,
Theorem \ref{2.9} implies that it vanishes.
\end{remark}

\section{Obstructions from almost flat relative classes} \label{sect3}

In this section we give obstructions to positive scalar 
curvature using almost flat classes in the relative setting.
Subsection \ref{subsect3.1} recalls the notions of relative cohomology and relative K-theory.  In Subsection \ref{subsect3.2} we discuss the
issue of basepoints for classifying maps in the relative setting, when the second space may not be
path connected.  Subsection \ref{subsect3.3} has the definition and basic properties of almost flat relative classes.
In Subsection \ref{subsect3.4} we use almost flat relative classes to give obstructions for a compact manifold-with-boundary
to admit a metric with positive scalar curvature and nonnegative mean curvature boundary.  Finally, in Subsection \ref{subsect3.5}
we give a local obstruction for a manifold to admit a metric with positive scalar curvature, using
almost flat relative classes.

\subsection{Relative cohomology} \label{subsect3.1}

Let $X$ be a compact topological space and let $Y$ be a closed subset of $X$.
Generators for the relative $K$-group $\KK^0(X,Y)$ are pairs $(E, {\sigma})$ where $E = E^+ \oplus E^-$ is a $\Z_2$-graded complex vector bundle on $X$ and ${\sigma} : E^+ \Big|_Y \rightarrow E^- \Big|_Y $ is an isomorphism
\cite[Chapter 2, Section 2.29]{Karoubi (1978)}. There is a Chern character
$\ch : \KK^0(X,Y) \rightarrow \HH^{even}(X, Y; \Q)$ \cite[Chapter 5, Section 3.26]{Karoubi (1978)}. We write $\ch(E, \sigma)$ for the Chern character
of the pair $(E, \sigma)$.

Now let $X$ be a compact topological space and let $Y$ be a 
compact topological space with a continuous map $\mu : Y \rightarrow X$.
Generators for the relative $K$-group $\KK^0(Y \rightarrow X)$ are pairs $(E, \sigma)$ where $E = E^+ \oplus E^-$ is a $\Z_2$-graded complex vector bundle on $X$ and $\sigma : \mu^* E^+ \rightarrow \mu^* E^- $ is an isomorphism.

Equivalently, in terms of mapping cylinders,
let $Z$ be the quotient of $X \cup \left( [0,1] \times Y \right)$ by the equivalence relation
$\mu(y) \sim (\{0\} \times y)$. Put $Z^\prime = \{1\} \times Y$, a subset of $Z$. Then 
$\KK^*(Y \rightarrow X) =
\KK^*(Z, Z^\prime)$. 

If $X^\prime$ and $Y^\prime$ are other such spaces, and we have 
maps $X \rightarrow X^\prime$ and $Y \rightarrow Y^\prime$ so that the diagram
\begin{equation} \label{3.1}
\begin{tikzcd}
    Y \arrow[r] \arrow{d} &  Y^\prime \arrow[d] \\
X \arrow[r] &  X^\prime
\end{tikzcd}
\end{equation}
commutes, we will denote it by a map
$\nu : (X, Y) \rightarrow (X^\prime, Y^\prime)$. There is a pullback
$\nu^* : \KK^*(Y^\prime \rightarrow X^\prime) \rightarrow
\KK^*(Y \rightarrow X)$.

As a variation, suppose that we additionally have a complex vector bundle
$V$ on $Y$, and an isomorphism ${\sigma} : \mu^* E^+  
\oplus V \rightarrow \mu^* E^- \oplus V$.  
We can find some $N \in \Z^+$ and a vector bundle
$G$ on $Y$ so that $V \oplus G$ 
can be identified with a trivial bundle $Y \times \C^N$. Put
$\widehat{E}^\pm = E^\pm \oplus (X \times \C^N)$. Define an isomorphism
$\widehat{\sigma}  : \mu^* \widehat{E}^+ \rightarrow \mu^* \widehat{E}^-$ 
by $\widehat{\sigma} = \sigma \oplus \Id \Big|_{G}$. We put
$[E, {V}, \sigma] = [\widehat{E}, \widehat{\sigma}]$ as an element of 
$\KK^0 \left( Y \rightarrow X \right)$; it is independent
of the choices made. We write $\ch(E, V, \sigma) = \ch(\widehat{E}, \widehat{\sigma})$.

We define the relative cohomology groups 
$\HH^*(Y \rightarrow X; \Q)$ by the
algebraic mapping cone construction.
That is, the relative cochains are
\begin{equation} \label{3.2}
C^k(Y \rightarrow X; \Q) =
C^k(X; \Q) 
 \oplus
C^{k-1}(Y; \Q)
\end{equation}
with differential 
$d(c^k, c^{k-1}) = (dc^k, \mu^* c^k - 
dc^{k-1})$.
 There are Chern characters
$\ch : \KK^0(Y \rightarrow X) \rightarrow 
\HH^{even}(Y \rightarrow X; \Q)$ and
$\ch : \KK^{-1}(Y \rightarrow X) \rightarrow 
\HH^{odd}(Y \rightarrow X; \Q)$.

\begin{remark} \label{3.3}
If $Y$ is a closed subset of
$X$ then $\KK^*(Y \rightarrow X) = \KK^*(X, Y)$ and
$\HH^*(Y \rightarrow X; \Q) = \HH^*(X, Y; \Q)$.
\end{remark}

\subsection{Basepoints} \label{subsect3.2}

When talking about classifying maps for pairs $(X, Y)$, where $Y$ may be disconnected, there is an
issue of how to deal with basepoints.  We first describe classifying maps in terms of fundamental
groupoids, where such issues do not arise.  We then show how to translate this in terms of 
fundamental groups.

If $h : \Gamma^\prime \rightarrow \Gamma$ is a homomorphism between discrete groups then we can
talk about
$\KK^*(B\Gamma^\prime \rightarrow B\Gamma)$ and
$\HH^*(B\Gamma^\prime \rightarrow B\Gamma; \Q)$.  In the special case when $h$ is injective,
there are two different algebraic notions of
the relative group cohomology \cite{Adamson (1954),Takasu (1959)}. Our $\HH^*(B\Gamma^\prime
\rightarrow B\Gamma; \Q)$ corresponds to the one defined in \cite{Takasu (1959)}.

If $\Gamma$ is replaced by a topological groupoid ${\mathcal G}$ then there is 
again a classifying space $B{\mathcal G}$
\cite{Segal (1968)}. Given a finite collection 
$\{ {\mathcal G}_\alpha \}$ of topological groupoids
with continuous homomorphisms $h_\alpha : {\mathcal G}_\alpha \rightarrow {\mathcal G}$
we can define $\KK^*(\coprod_\alpha B {\mathcal G}_\alpha \rightarrow B{\mathcal G})$.

The groupoids ${\mathcal G}$ that are relevant for us are transitive \'etale groupoids.
Here ``transitive'' means that the range and source maps $r, s : {\mathcal G}^{(1)}
\rightarrow {\mathcal G}^{(0)}$ are such that $(r,s) : 
{\mathcal G}^{(1)}
\rightarrow {\mathcal G}^{(0)} \times {\mathcal G}^{(0)}$ is surjective. Then ${\mathcal G}$ is Morita equivalent to a discrete
group, namely the isotropy group ${\mathcal G}_p = r^{-1}(p) \cap s^{-1}(p)$ for any
$p \in {\mathcal G}^{(0)}$.

We now describe the classifying map $\nu$ for pairs, arising from fundamental groupoids.
Let $X$ be a compact path connected topological space.  Let
${\mathcal G}_X$ be the fundamental groupoid of $X$ \cite[Chapter 2.5]{May (1999)} and let
$B{\mathcal G}_X$ be its classifying space, as constructed in the sense of classifying spaces of topological categories
\cite{Segal (1968)}.
There is an inclusion $X \rightarrow B{\mathcal G}_X$.
Let $Y$ be a closed subset of $X$, with path components
$\{Y_\alpha \}$.
We obtain
a commutative diagram
\begin{equation} \label{3.4}
\begin{tikzcd}
Y  \arrow[r] \arrow[d] &   \coprod_\alpha B{\mathcal G}_{Y_\alpha} \arrow[d] \\
X \arrow[r] &  B{\mathcal G}_X. 
\end{tikzcd}
\end{equation}
For brevity, we will write such a diagram as a map
$(X, Y) \rightarrow \left( B{\mathcal G}_X, \coprod_\alpha B{\mathcal G}_{Y_\alpha}
\right)$ and will call it the classifying map for the pair $(X, Y)$.

To write the classifying map in terms of groups, instead of groupoids,
choose a basepoint $x_0 \in X$ and basepoints $y_\alpha \in 
Y_\alpha$. Let $\gamma_\alpha$ be a path from
$x_0$ to $y_\alpha$. This gives a homomorphism 
$\pi_1(Y_\alpha, y_\alpha) \rightarrow \pi_1(X, x_0)$.  The path groupoid ${\mathcal G}_X$ is Morita-equivalent to $\pi_1(X, x_0)$, and similarly for ${\mathcal G}_{Y_\alpha}$. We obtain a diagram
\begin{equation} \label{3.5}
\begin{tikzcd}
    \coprod_\alpha B{\mathcal G}_{Y_\alpha} \arrow[r] \arrow{d} &  \coprod_\alpha B\pi_1(Y_\alpha, y_\alpha) \arrow[d] \\
B{\mathcal G}_X \arrow[r] &  B\pi_1(X, x_0),
\end{tikzcd}
\end{equation}
depending on $\{ \gamma_\alpha \}$, that commutes up to homotopy.

Putting (\ref{3.4}) and (\ref{3.5}) together, we obtain a diagram
\begin{equation} \label{3.6}
\begin{tikzcd}
    Y \arrow[r] \arrow{d} &  \coprod_\alpha B\pi_1(Y_\alpha, y_\alpha) \arrow[d] \\
X \arrow[r] &  B\pi_1(X, x_0),
\end{tikzcd}
\end{equation}
depending on $\{ \gamma_\alpha \}$, that commutes up to homotopy. 
For brevity, we will write such a diagram as a map
$(X, Y) \rightarrow \left( B\pi_1(X, x_0), \coprod_\alpha B\pi_1(Y_\alpha, y_\alpha)
\right)$.
Different choices of $\{ \gamma_\alpha \}$ will give different, but equivalent,
groups and diagrams.

Let $\Gamma$ be a discrete group and let $\{ \Gamma^\prime_\alpha \}$ be a finite collection of discrete
groups with homomorphisms $h_\alpha : \Gamma^\prime_\alpha \rightarrow \Gamma$. Given 
homomorphisms $\pi_1(X, x_0) \rightarrow \Gamma$ and
$\pi_1(Y_\alpha, y_\alpha) \rightarrow \Gamma^\prime_\alpha$
so that the diagrams
\begin{equation} \label{3.7}
\begin{tikzcd}
    \pi_1(Y_\alpha, y_\alpha) \arrow[r] \arrow{d} &  \Gamma^\prime_\alpha \arrow[d] \\
\pi_1(X, x_0) \arrow[r] &  \Gamma
\end{tikzcd}
\end{equation}
commute, we
obtain a diagram
\begin{equation} \label{3.8}
\begin{tikzcd}
    Y \arrow[r] \arrow{d} &  \coprod_\alpha B\Gamma^\prime_\alpha \arrow[d] \\
X \arrow[r] &  B\Gamma,
\end{tikzcd}
\end{equation}
which we will write as a map
$\nu : (X, Y) \rightarrow
(B\Gamma,  B \Gamma^\prime)$. It is defined up to homotopy.
Note that the maps $B \Gamma^\prime_\alpha \rightarrow B\Gamma$ need not be inclusions,
so $\nu$ is not a map of pairs in the conventional sense.  Rather, it is to 
be interpreted as the diagram (\ref{3.8}). There is still a pullback map
$\nu^* : \KK^*( \coprod_\alpha B\Gamma^\prime_\alpha \rightarrow B\Gamma) \rightarrow
\KK^*( Y \rightarrow X) \cong \KK^*(X, Y)$.

\subsection{Almost flat relative classes} \label{subsect3.3}

\begin{definition} \label{3.9}
Let ${\mathcal N}$ be a compact manifold-with-boundary.  Let ${\mathcal N}^\prime$
be a compact manifold-with-boundary with a smooth map
$\mu : {\mathcal N}^\prime \rightarrow {\mathcal N}$.
Put Riemannian metrics on ${\mathcal N}$ and ${\mathcal N}^\prime$.
Given $\beta \in \KK^0({\mathcal N}^\prime\rightarrow {\mathcal N})$, 
we say that $\beta$ is almost flat if for every
$\epsilon > 0$,
we can find
\begin{itemize}
    \item A $\Z_2$-graded Hermitian vector bundle $E$ on ${\mathcal N}$
    and an isometric isomorphism $\sigma : \mu^* E^+ 
    \rightarrow \mu^* E^-$ with
    $(E, \sigma )$ representing $\beta$, and
    \item A Hermitian connection $\nabla^E$ on $E$ whose curvature satisfies
    $\| F^\pm \| \le \epsilon$, and also $\|(\mu^* \nabla^E) \sigma \| \le \epsilon$.
\end{itemize}
\end{definition}

Let $\KK^0_{af}({\mathcal N}^\prime \rightarrow {\mathcal N})$ denote the almost flat classes.
The above definition of an almost flat relative class is essentially the same as that in
\cite[Definition 3.17]{Kubota (2022)}. The definition in \cite{Kubota (2022)} is effectively when 
${\mathcal N}^\prime$ is a submanifold of ${\mathcal N}$.

\begin{definition} \label{3.10}
Let $\Gamma$ be a discrete group and let $\{\Gamma^\prime_\alpha\}$ be a finite collection of
discrete groups, with homomorphisms $h_\alpha : \Gamma^\prime_\alpha \rightarrow \Gamma$.
Given $\eta \in \KK^{0} ( \coprod_\alpha B\Gamma^\prime_\alpha \rightarrow B\Gamma)$,
we say that $\eta$ is almost flat if 
for every choice of compact manifold-with-boundary ${\mathcal N}$, compact manifolds-with-boundary
$\{{\mathcal N}^\prime_\alpha \}$, smooth maps
$\mu_\alpha : {\mathcal N}^\prime_\alpha \rightarrow {\mathcal N}$
and continuous map of pairs
$\nu : ({\mathcal N}, {\mathcal N}^\prime) \rightarrow (B\Gamma,  B\Gamma^\prime)$,
we have $\nu^* \eta \in \KK^0_{af}(\coprod_\alpha  {\mathcal N}^\prime_\alpha \rightarrow {\mathcal N})$.
\end{definition}

Let $\KK^{0}_{af}(\coprod_\alpha B\Gamma^\prime_\alpha \rightarrow B\Gamma)$ denote the almost flat classes. Note that if $B\Gamma$ and $\{B\Gamma^\prime_\alpha\}$ are smooth compact manifolds then a class
$\eta \in \KK^{0}(\coprod_\alpha B\Gamma^\prime_\alpha \rightarrow B\Gamma)$ is almost flat in
the sense of Definition \ref{3.10} if and
only if it is almost flat in the sense of Definition \ref{3.9}.

For the odd case, define the almost flat elements
$\KK^{-1}_{af}(\coprod_\alpha B\Gamma^\prime_\alpha \rightarrow B\Gamma)$ to be the elements of
$\KK^{-1}(\coprod_\alpha B\Gamma^\prime_\alpha \rightarrow B\Gamma)$
whose product with
$[B\Z] \in \KK^{-1}(B\Z)$ lies in 
$\KK^{0}_{af}(\coprod_\alpha B(\Gamma^\prime_\alpha \times \Z) \rightarrow 
B(\Gamma \times \Z))$.

\begin{example} \label{3.11}
Let $\Gamma$ be a discrete group and let 
$\{\Gamma^\prime_\alpha\}$ be a finite collection of
trivial groups. The exact sequence in $\KK$-theory implies that
\begin{equation} \label{3.12}
\KK^{0}(\coprod_\alpha B\Gamma^\prime_\alpha \rightarrow B\Gamma) = 
\Ker \left( \KK^0(B\Gamma) \rightarrow \oplus_\alpha \KK^0(B\Gamma^\prime_\alpha) \right).
\end{equation}
Suppose that $\KK^{0}_{af}(B\Gamma) =
\KK^{0}(B\Gamma)
$. As
$B\Gamma^\prime_\alpha$ maps to a point in $B\Gamma$, 
it follows that
$\KK^{0}_{af}(\coprod_\alpha B\Gamma^\prime_\alpha \rightarrow B\Gamma) =
\KK^{0}(\coprod_\alpha B\Gamma^\prime_\alpha \rightarrow B\Gamma)
$.
\end{example}

\begin{example} \label{3.13}
Suppose that $\Gamma = \{e\}$ and $\Gamma_1 = \Z$. Then $B\Gamma = \pt$ and $B\Gamma_1 = S^1$.
From the exact sequence in $\KK$-theory, the boundary map $\KK^{-1}(S^1) \rightarrow
\KK^0(S^1 \rightarrow \pt)$ is an isomorphism, so 
$\KK^0(S^1 \rightarrow \pt) \cong \Z$. To analyze
$\KK^0_{af}(S^1 \rightarrow \pt)$,
we take ${\mathcal N} = \pt$ and
${\mathcal N}^\prime_1 = S^1$. The vector bundle $E$ is just a $\Z_2$-graded inner product space on 
$\pt$ and $\nabla^E$ is trivial.  The pullback $\mu_1^* E$ is a trivial $\Z_2$-graded
Hermitian vector bundle on $S^1$ and $\mu_1^* \nabla^E$ is the trivial connection.  If
$\phi_1 : \mu_1^* E^+ \rightarrow \mu_1^* E^-$ is an isometric isomorphism then we can
identify $\phi_1$ with a map $\Phi : S^1 \rightarrow U(N)$ satisfying
$\|\Phi^{-1} \Phi^\prime\| \le \epsilon$.  If $\epsilon$ is small, independent of $N$, then this forces $\det \Phi$ to have vanishing winding
number (one can reduce to when $\Phi$ takes value in $U(1)^N$), so $\KK^0_{af}(S^1 \rightarrow \pt) = 0$.
\end{example}

\begin{remark} \label{3.14}
Example \ref{3.13} shows that if $\Gamma^\prime \rightarrow \Gamma$ is not injective then we cannot
expect that $\KK^{*}_{af}(\coprod_\alpha B\Gamma^\prime_\alpha \rightarrow B\Gamma) \otimes \Q = 
\KK^{*}(\coprod_\alpha B\Gamma^\prime_\alpha \rightarrow B\Gamma) \otimes \Q$. Put
$\widetilde{\Gamma}^\prime_\alpha = {\Gamma}^\prime_\alpha/\Ker (h_\alpha)$. The induced 
homomorphism $\widetilde{h}_\alpha : \widetilde{\Gamma}^\prime_\alpha \rightarrow \Gamma$
is injective.  As there is a map $\KK^*_{af}
(\coprod_\alpha B \widetilde{\Gamma}^\prime_\alpha \rightarrow B\Gamma) \rightarrow
\KK^*_{af}
(\coprod_\alpha B\Gamma^\prime_\alpha \rightarrow B\Gamma)$, by switching
from $\Gamma^\prime$ to $\widetilde{\Gamma}^\prime$ we can effectively
just consider almost flat relative $\KK$-theory classes that come from $h_\alpha$ being injective.
It seems conceivable that the map $\KK^{*}_{af}(\coprod_\alpha B\Gamma^\prime_\alpha \rightarrow B\Gamma) \otimes \Q 
\rightarrow
\KK^{*}(\coprod_\alpha B\Gamma^\prime_\alpha \rightarrow B\Gamma) \otimes \Q$ is an isomorphism if the $h_\alpha$'s are
injective, at least when $B\Gamma$ and the $B \Gamma^\prime_\alpha$'s are finite CW-complexes.
\end{remark}

\subsection{Manifolds with boundary} \label{subsect3.4}

\begin{theorem} \label{3.15}
Let $M$ be a compact Riemannian spin manifold-with-boundary. 
If $M$ has positive scalar curvature and $\partial M$ has nonnegative mean curvature
then for all $\beta \in \KK^*_{af} 
\left(M,  \partial M \right)$,
we have
\begin{equation} \label{3.16}
\int_{M} \widehat{A}(TM) \wedge 
\ch (\beta) = 0.
\end{equation}
\end{theorem}
\begin{proof}
Suppose first that $M$ is even dimensional.  Let $S$ denote the $\Z_2$-graded
spinor bundle on $M$. As in Definition \ref{3.9}, let $E$ be a $\Z_2$-graded Hermitian vector bundle
on $M$, with an isometric isomorphism $\sigma : E^+ \big|_{\partial M} \rightarrow
E^- \big|_{\partial M}$,
so that $(E, \sigma)$ represents $\beta$. Given $\epsilon > 0$, let $\nabla^E$ be a Hermitian connection on
$E$ so that $\|F^{\pm} \| \le \epsilon$ and $\| \nabla^E \sigma \| \le \epsilon$. 
Give $S \otimes E$ the total $\Z_2$-grading.
Let $D^E$
be the Dirac-type operator on $C^\infty(M; S \otimes E)$. Then
\begin{equation} \label{3.17}
(D^E)^2 = \nabla^* \nabla + \frac14 R - \frac14 \left[ \gamma^\mu, \gamma^\nu \right] 
F^E_{\mu \nu}. 
\end{equation}

Define a projection $\Pi$ on $C^\infty \left( \partial M; (S \otimes E) \big|_{\partial M} \right)$
by 
\begin{equation} \label{3.18}
\Pi = \begin{pmatrix}
0 & \sqrt{-1} \epsilon_S \gamma^n \sigma^{-1} \\
\sqrt{-1} \epsilon_S \gamma^n \sigma & 0
\end{pmatrix}.
\end{equation}
Let ${\mathcal D}$ be the self-adjoint Dirac-type operator, densely defined on
$L^2(M; S \otimes E)$, with boundary condition $(\Pi - I) \psi \big|_{\partial M} = 0$.
As in the proof of Lemma \ref{2.19}, if $\epsilon$ is small enough then
$\Ker({\mathcal D}) = 0$.

To compute the index of ${\mathcal D}$, we use a doubling argument as in the proof of
Lemma \ref{2.35}.  Because of the local boundary conditions, the index only depends on the principal
symbols of $D^E$ and $\Pi$, up to homotopy. We deform the Riemannian metric and the Hermitian connection $\nabla^E$
so that they are products near $\partial M$, and we can also deform so that
$\nabla^{E^-} \big|_{\partial M} = \sigma \circ \nabla^{E^+} \big|_{\partial M} \circ \sigma^{-1}$.
The form $\tr_s e^{\frac{i}{2 \pi} F^E}$ vanishes in a neighborhood of $\partial M$ and represents
$\ch(\beta)$.

We pass to the double $DM$, whose $\Z_2$-involution has fixed point set $\partial M$. We can extend $S \otimes E$ to $DM$ so that the generator of the
$\Z_2$-involution acts on $(S \otimes E) \big|_{\partial M}$ by $i \epsilon_S \gamma^n \otimes
\begin{pmatrix}
0 & \sigma^{-1} \\
\sigma & 0
\end{pmatrix}
$. The index of ${\mathcal D}$ is the same as the $\Z_2$-invariant index on $DM$, which is
the same as the orbifold index on the orbifold $DM/\Z_2$. From the orbifold index theorem, the index equals 
$\int_M \widehat{A}(TM) \wedge \tr_s e^{\frac{i}{2 \pi} F^E} =
\int_M \widehat{A}(TM) \wedge \ch(\beta)$.

This proves the theorem in the even dimensional case.  In the odd dimensional case, we take a
product with $S^1$ to reduce to the even dimensional case.
\end{proof}

\begin{corollary} \label{3.19}
Let $M$ be a compact Riemannian spin manifold-with-boundary. Let $\{ \partial_\alpha M \}$ be
the boundary components. Suppose that $M$ has positive scalar curvature and
$\partial M$ has nonnegative mean curvature.

Let $\Gamma$ be a discrete group and let $\{\Gamma^\prime_\alpha\}$ be a finite collection of
discrete groups, with homomorphisms $h_\alpha : \Gamma^\prime_\alpha \rightarrow \Gamma$.
Given $\eta \in \KK^*_{af} \left( \coprod_\alpha B\Gamma^\prime_\alpha  \rightarrow B\Gamma \right)$ and a continuous map of pairs
$\nu :  
(M, \partial M) \rightarrow (B\Gamma, B\Gamma^\prime)$,
we have
\begin{equation} \label{3.20}
\int_{M} \widehat{A}(TM) \wedge 
\ch(\nu^* \eta) = 0.
\end{equation}
\end{corollary}

\begin{remark}
There is some intersection between Corollary \ref{3.19} and \cite[Theorem 2.19]{Baer-Hanke (2023)}
when each $\Gamma^\prime_\alpha$ is trivial. The proof in \cite[Theorem 2.19]{Baer-Hanke (2023)} uses APS boundary
conditions.
\end{remark}

\begin{corollary} \label{3.21}
Let $M$ be a compact Riemannian spin manifold-with-boundary. Let $\{ \partial_\alpha M \}$ be
the boundary components. Suppose that $M$ has positive scalar curvature and
$\partial M$ has nonnegative mean curvature.

Put $\Gamma = \pi_1(M)$ and $\Gamma^\prime_\alpha = \pi_1(\partial_\alpha M)$, where basepoints are
handled as in Section \ref{subsect3.2}. Put $\widetilde{\Gamma}^\prime_\alpha = \Gamma^\prime_\alpha/\Ker( \Gamma^\prime_\alpha
\rightarrow \Gamma)$. 
Let $\nu$ be the composition
$(M, \partial M) \rightarrow (B\Gamma, B{\Gamma}^\prime) \rightarrow (B\Gamma, B\widetilde{\Gamma}^\prime)$.
Suppose that $\KK^*_{af} \left( \coprod_\alpha B \widetilde{\Gamma}^\prime_\alpha  \rightarrow B\Gamma \right) \otimes \Q = \KK^* \left( \coprod_\alpha B \widetilde{\Gamma}^\prime_\alpha  \rightarrow B\Gamma \right) \otimes \Q$.
Then $\nu_* [M, \partial M]$ vanishes in $\HH_{\dim(M)}(B\Gamma, B\widetilde{\Gamma}^\prime; \Q)$.
\end{corollary}
\begin{proof}
Given Corollary \ref{3.19}, the proof is similar to that of Corollary \ref{2.7}.
    \end{proof}

    \begin{remark} \label{3.22}
    The motivation to pass from $\Gamma^\prime$ to $\widetilde{\Gamma}^\prime$ comes from Remark \ref{3.14}.
    \end{remark}

\begin{remark} \label{3.23}
To put Corollary \ref{3.21} in perspective,
we discuss the index theory obstructions to the existence of positive scalar curvature metrics on compact 
manifolds-with-boundary.

Let $M$ be a compact connected spin manifold-with-boundary.  For simplicity, we just discuss the
case when $\partial M$ is nonempty and connected; the issue of basepoints in the general case can be handled 
as in Section \ref{subsect3.2}. Choose a basepoint $m_0 \in \partial M$, put $\Gamma = \pi_1(M, m_0)$
and put $\Gamma^\prime = \pi_1(\partial M, m_0)$.

Suppose first that we are given a
metric $g_{\partial M}$ on $\partial M$ with $R_{\partial M} > 0$, and we want to know if we can extend it to a metric on $M$ 
so that $R_M > 0$ and $H_{\partial M} \ge 0$. Choose an extension to an arbitrary metric
$g_M$ on $M$. Since $R_{\partial M} > 0$, 
the Dirac operator on $M$, with APS boundary conditions, has a well defined index in 
$\KK_*(C^*_{max} \Gamma)$, the $K$-theory of the maximal group $C^*$-algebra. 
If $g_M$ has 
$R_M > 0$ and $H_{\partial M} \ge 0$ then the index vanishes. 

If we ask for the obstruction to finding {\em some} metric $g_M$ on $M$ with $R_M > 0$, 
$H_{\partial M} \ge 0$ and $R_{\partial M} > 0$, then we are effectively allowing ourselves to
perturb a given boundary metric $g_{\partial M}$ with $R_{\partial M} > 0$ to another boundary
metric with positive scalar curvature.  This changes the index by an element of
$\Image \left( \KK_*(C^*_{max} \Gamma^\prime) \rightarrow \KK_*(C^*_{max} \Gamma)\right)$.
Hence the obstruction lies in $\Coker \left( \KK_*(C^*_{max} \Gamma^\prime) \rightarrow \KK_*(C^*_{max} \Gamma)\right)$.

There is a relative group $C^*$-algebra $C^*_{max}(\Gamma, \Gamma^\prime)$ that fits into a
short exact sequence
\begin{equation} \label{3.24}
0 \rightarrow S C^*_{max} \Gamma \rightarrow C^*_{max}(\Gamma, \Gamma^\prime) \rightarrow C^*_{max} \Gamma^\prime \rightarrow 0.
\end{equation}
From the long exact sequence of $\KK$-groups,
we can also think of the obstruction as lying in $\Ker \left( \KK_{*-1}(C^*_{max} 
(\Gamma, \Gamma^\prime)) \rightarrow \KK_{*-1}(C^*_{max} \Gamma^\prime)\right)$. This obstruction was constructed, as an element of
$\KK_{*-1}(C^*_{max} (\Gamma, \Gamma^\prime))$, in \cite{Chang-Weinberger-Yu (2020)}.  The equivalence with the
above use of the APS boundary conditions was shown in \cite{Schick-Seyedhosseini (2021)}. The fact that the element maps to zero in $\KK_{*-1}(C^*_{max} \Gamma^\prime)$ is a 
reflection of the assumption that $R_{\partial M} > 0$.

It is known that the existence of a metric $g_M$ with $R_M > 0$, $H_{\partial M} \ge 0$ and
$R_{\partial M} > 0$ is equivalent to the existence of such a metric which is a product near
$\partial M$ \cite{Baer-Hanke (2023)}. In fact, the papers 
\cite{Chang-Weinberger-Yu (2020),Schick-Seyedhosseini (2021)} consider metrics that are products near
$\partial M$.

We now remove the assumption that $R_{\partial M} > 0$. One obstruction to having $R_M > 0$ and
$H_{\partial M} \ge 0$ comes from Theorem \ref{3.15}.
Another approach is to take the double $DM$ and smooth the resulting metric to get a 
$\Z_2$-invariant Riemannian
metric on $DM$ with $R_{DM} > 0$ \cite{Rosenberg-Weinberger}. Now $\pi_1(DM, m_0)$ is the
amalgamated free product
$\pi_1(M, m_0) \star_{\pi_1(\partial M, m_0)} \pi_1(M, m_0)$.  Since 
the homomorphism $\pi_1(\partial M, m_0)
\rightarrow \pi_1(M, m_0)$ may not be injective, put
$\Gamma = \pi_1(M, m_0)$ and $\widetilde{\Gamma}^\prime = 
\pi_1(\partial M, m_0)/\Ker(\pi_1(\partial M, m_0)
\rightarrow \pi_1(M, m_0))$. Then the induced homomorphism $\widetilde{\Gamma}^\prime \rightarrow \Gamma$ is injective and $\pi_1(DM, m_0) = \Gamma \star_{\widetilde{\Gamma}^\prime} \Gamma$. 

There are $\Z_2$-equivariant maps $\KK_*(DM) \rightarrow \KK_*(B\Gamma \cup_{B\widetilde{\Gamma}^\prime} B\Gamma)) \rightarrow \KK_*(C^*_{max}(\Gamma \star_{\widetilde{\Gamma}^\prime} \Gamma))$. Taking into account that the involution on $DM$ is orientation reversing, we can pass to the $\Z_2$-quotient to get maps
\begin{equation} \label{3.25}
\KK_*(M, \partial M) \rightarrow \KK_*(B\Gamma, B\widetilde{\Gamma}^\prime) \rightarrow
\KK_*^-(C^*_{max}(\Gamma \star_{\widetilde{\Gamma}^\prime} \Gamma)),
\end{equation}
where
$\KK_*^-(C^*_{max}(\Gamma \star_{\widetilde{\Gamma}^\prime} \Gamma))$ denotes the part of
$\KK_*(C^*_{max}(\Gamma \star_{\widetilde{\Gamma}^\prime} \Gamma))$ that is odd under involution.
Under the curvature assumption, the strong Novikov conjecture
for $C^*_{max}$ (formulated $\Z_2$-equivariantly)
implies that the image of the fundamental class $[M, \partial M] \in \KK_*(M, \partial M)$ in
$\KK_*(B\Gamma, B\widetilde{\Gamma}^\prime) \otimes \Q$ vanishes; compare with Corollary \ref{3.21}.
If $\widetilde{\Gamma}^\prime = \Gamma$ then it is known that a metric with $R_M > 0$ and $H_{\partial M} \ge 0$
always exists \cite[Theorem 2.2]{Rosenberg-Weinberger}.
\end{remark}

\subsection{Local obstruction using almost flat relative bundles} \label{subsect3.5}

In this subsection we use almost flat relative bundles to give localized obstructions to
having positive scalar curvature.  We make the following assumption.
\begin{assumption} \label{3.26}
Given $r_0, D > 0$, put $r_0^\prime = \frac{1}{256} r_0^2 D^2$
and $D^\prime = D + \frac{32}{r_0 D}$. 
Let $M$ be a Riemannian spin 
manifold, possibly incomplete. Let $K$ be a compact submanifold-with-boundary of $M$.
Suppose that 
\begin{itemize}
\item The distance neighborhood $N_{D^\prime}(K)$ lies in a compact 
submanifold-with-boundary ${\mathcal C}$ of $M$,
\item $R > 0$ on $K$,
\item $R \ge r_0$ on $N_{D}(K) - K$ and
\item $R \ge - r_0^\prime$ on $N_{D^\prime}(K) - N_{D}(K)$.
\end{itemize}
\end{assumption}

\begin{theorem} \label{3.27}
Suppose that Assumption \ref{3.26} holds.
We can identify $\KK^* \left( ({\mathcal C} - \Int(K))
\rightarrow {\mathcal C} \right)$ with
$\KK^*\left({\mathcal C},  {\mathcal C} - \Int(K) \right)$.
Then given $\beta \in \KK^*_{af} 
\left({\mathcal C},  {\mathcal C} - \Int(K)
\right)$,
we have
\begin{equation} \label{3.28}
\int_{\mathcal C} \widehat{A}(TM) \wedge 
\ch (\beta) = 0.
\end{equation}
\end{theorem}
\begin{proof}
Suppose first that $M$ is even dimensional.  Let $\epsilon > 0$ be a
small parameter, which we will adjust.  Let ${\mathcal N}^\epsilon$
and $f$ be as in the proof of Theorem \ref{2.9}. (We now have $\partial M = \emptyset$.)
As in Definition \ref{3.9}, let $(E, \sigma)$ be a pair that represents 
$\beta \big|_{{\mathcal N}^\epsilon}$,
where $E$ is a $\Z_2$-graded Hermitian vector bundle on ${\mathcal N^\epsilon}$
and $\sigma : E^+ \big|_{{\mathcal N}^\epsilon - \Int(K)} \rightarrow
E^- \big|_{{\mathcal N}^\epsilon - \Int(K)}$ is an isometric isomorphism.

Put
\begin{equation} \label{3.29}
A = 
\begin{pmatrix}
\nabla^{+} & f \sigma^{-1} \\
f\sigma & \nabla^{-}
\end{pmatrix},
\end{equation}
a superconnection on $E$. 
Let $D^E$ be the quantization of $A$, i.e.
\begin{equation} \label{3.30}
D^E = 
\begin{pmatrix}
D^+ & \epsilon_S f \sigma^{-1} \\
\epsilon_S f \sigma & D^-
\end{pmatrix},
\end{equation}
where $D^\pm$ is the Dirac operator on $C^\infty({\mathcal N}^\epsilon; S \otimes E^\pm)$ and $\epsilon_S$ is the $\Z_2$-grading
operator on spinors on ${\mathcal N}^\epsilon$. Then
\begin{align} \label{3.31}
(D^E)^2 = 
& \begin{pmatrix}
(D^+)^2 + f^2 & 0 \\
0 & (D^-)^2 + f^2
\end{pmatrix} + 
\sqrt{-1} \epsilon_S c(df)
\begin{pmatrix}
0 &  \sigma^{-1} \\
 \sigma & 0
\end{pmatrix} + \\
& \sqrt{-1} \epsilon_S f
\begin{pmatrix}
0 & c(\nabla^{+} \circ \sigma^{-1} - \sigma^{-1} \circ \nabla^{-}) \\
 c(\nabla^{-} \circ \sigma - \sigma \circ \nabla^{+}) & 0
\end{pmatrix}. \notag
\end{align}

We impose the boundary condition $\Pi \psi = \psi$ on $\partial {\mathcal N}^\epsilon$,
as in the proof of
Theorem \ref{2.9}.  Let ${\mathcal D}$ denote the ensuing self-adjoint
operator.
As in Lemma \ref{2.19}, 
there is an $\epsilon_0 > 0$ so that if $\epsilon < \epsilon_0$, and in addition
$\| F^\pm \| \le \epsilon$ and $\| \nabla^E \sigma \| \le \epsilon$, then
the kernel of ${\mathcal D}$ vanishes.

To compute the index of ${\mathcal D}$, we follow the method of proof of
Lemma \ref{2.35}. We deform the metric on ${\mathcal N}^\epsilon$ so that
it is a product near the boundary and we deform $f$ to zero.
We deform $\nabla^+$ to be a product near the boundary.  We deform
$\nabla^-$ to be $\sigma \circ \nabla^+ \circ \sigma^{-1}$ near
the boundary.  By the doubling argument, the index of ${\mathcal D}$
equals $\int_{{\mathcal N}^\epsilon} \widehat{A}(T{\mathcal N}^\epsilon)
\wedge \tr_s \left( e^{\frac{i}{2\pi} (\nabla^E)^2} \right)$.
The form $\tr_s \left( e^{\frac{i}{2\pi} (\nabla^E)^2} \right)$ vanishes
near $\partial {\mathcal N}^\epsilon$ and represents 
$\ch \left( \beta \big|_{{\mathcal N}^\epsilon} \right)$. Extending it by zero to ${\mathcal C}$ gives a
representative of $\ch(\beta)$; note that ${\mathcal N} - {\mathcal N}^\epsilon$
is contained in ${\mathcal N} - \Int(K)$.
Hence the index equals
\begin{equation} \label{3.32}
\int_{{\mathcal N}^\epsilon} \widehat{A}(T{\mathcal N}^\epsilon)
\wedge \tr_s \left( e^{\frac{i}{2\pi} (\nabla^E)^2} \right) = \int_{{\mathcal C}} \widehat{A}(TM)
\wedge \ch(\beta).
\end{equation}
This proves the theorem if $M$ is even dimensional.  If $M$ is odd dimensional, consider
$M^\prime = M \times S^1$
and $K^\prime = 
K \times S^1$. Then we reduce to the even dimensional case. 
\end{proof}

\begin{example} \label{3.34}
Let $Y$ be a compact connected spin manifold of even dimension. 
Let $F$ be a finite subset of $Y$.
Let $M = Y \#_F Z$ be the result of taking the connect sum of $Y$ with a
possibly noncompact and possibly disconnected spin manifold $Z$, over the points of $F$.
(If $F$ is a point then this would be the usual connect sum.)
We assume that the connect sum is taken within a small neighborhood $U$ of $F$.
Let $\tau : M \rightarrow Y$ be a map which is the identity on $Y-U \subset M$ and maps 
$M-(Y-U)$ to $U$, while sending the complement of a compact subset $K \subset M$ to $F$.

Given an element $\delta \in \Ker \left( \KK^0(Y) \rightarrow \KK^0(F) \right)$,
suppose that $\delta$ lies in $\KK^0_{af}(Y)$.  Representing $\delta$ by a $\Z_2$-graded Hermitian
vector bundle $\widehat{E}$ with connection $\nabla^{\widehat{E}}$, put $E = 
\tau^* \widehat{E}$ and $\nabla^E = \tau^* \nabla^{\widehat{E}}$.
There is an isometric isomorphism $\sigma : E^+ \big|_{M - \Int(K)} \rightarrow
E^- \big|_{M - \Int(K)}$, pulled back from $F$, with $\nabla^E \sigma = 0$.
Suppose that $M$ has a Riemannian metric which satisfies Assumption \ref{3.26}. We have constructed an element of
$\KK^0_{af} 
\left({\mathcal C},  {\mathcal C} - \Int(K)
\right)$ and Theorem \ref{3.27} implies that
$\int_Y \widehat{A}(TY) \wedge \ch(\delta) = 0$. If $Y$ is odd dimensional 
and $\delta \in \KK^{-1}_{af}(Y)$ then we can
take the product with $S^1$ to obtain the same result.

In particular, if $M$ is complete with positive scalar curvature then by choosing
appropriate ${\mathcal C}$, we conclude that $\int_Y \widehat{A}(TY) \wedge \ch(\delta) = 0$.
This can be compared with
\cite[Theorem 0.2]{Wang-Zhang (2022)}.
\end{example}

\begin{example}
We show that Theorem \ref{3.27}, or more precisely its proof, gives a localized version of
\cite[Theorem 6.12]{Gromov-Lawson (1983)}.
Suppose first that $\dim(M)$ is even.
With reference to Assumption \ref{3.26}, given a finite cover 
$p : \widehat{\mathcal C} \rightarrow {\mathcal C}$ of ${\mathcal C}$, put
$\widehat{K} = p^{-1}(K)$. 
Suppose that for every $\epsilon > 0$, there is a  finite cover 
$p : \widehat{\mathcal C} \rightarrow {\mathcal C}$
and a map 
$m : \Int(\widehat{K}) \rightarrow S^n$ so that
\begin{itemize}
\item $m$ sends the complement of some compact subset of 
$\Int(\widehat{K})$ to a point,
\item $m$ is $\epsilon$-contracting on $2$-forms, and
\item $m$ has a nonzero degree $d$.
\end{itemize}
Extend $m$ to be a point map on $\widehat{\mathcal C} - \Int(\widehat{K})$.
Let $W$ be a $\Z_2$-graded Hermitian vector bundle on $S^n$ 
with a fixed Hermitian connection, so that
$\ch(W) = e [S^n] \in \HH^n(S^n; \Q)$ for some $e \in \Q - \{0\}$.
Define a $\Z_2$-graded Hermitian vector bundle
on ${\mathcal C}$ by $E = p_* m^* W$, i.e. for $x \in {\mathcal C}$, the fiber
$E_x$ is $\oplus_{y \in p^{-1}(x)} (m^*W)_y$.
Then
the curvature of the induced connection on $E$ is $O(\epsilon)$
in magnitude, and it vanishes outside of $\Int(K)$

Note that
although the connection on $E$ is flat on ${\mathcal C} - \Int(K)$, it will generally
have nontrivial holonomy there. As $m$ sends $\widehat{\mathcal C} - \Int(\widehat{K})$ to a
point, and $W^+ \big|_{\pt} \cong W^- \big|_{\pt}$,
there is a parallel isometric isomorphism
$\sigma : E^+ \big|_{{\mathcal C} - \Int(K)} \rightarrow E^- \big|_{{\mathcal C} - \Int(K)}$.
Then $\ch(E, \sigma) = d e [{\mathcal C}, {\mathcal C} - \Int(K)] \in 
\HH^n({\mathcal C}, {\mathcal C} - \Int(K); \Q)$.

If Assumption \ref{3.26} holds then the proof of Theorem \ref{3.27}
implies that
$\int_{\mathcal C} \widehat{A}(TM) \wedge \ch(E, \sigma) = 0$. However,
\begin{equation}
    \int_{\mathcal C} \widehat{A}(TM) \wedge \ch(E, \sigma) = 
\int_{\mathcal C} d e [{\mathcal C}, {\mathcal C} - \Int(K)] = de \neq 0,
\end{equation}
which is a contradiction.
If $\dim(M)$ is odd then we can take the product with a circle to get the same conclusion.

Hence the scalar curvature assumption cannot be satisfied.
This can be viewed as a localized version of \cite[Theorem 6.12]{Gromov-Lawson (1983)},
with the difference that we only consider finite covers in defining
$\Lambda^2$-enlargeability, as in \cite{Gromov-Lawson (1980)}.
\end{example}

\begin{corollary} \label{3.36}
Suppose that Assumption \ref{3.26} holds.
Let $\Gamma$ be a discrete group and let $\{\Gamma^\prime_\alpha\}$ be a finite collection of
discrete groups, with homomorphisms $h_\alpha : \Gamma^\prime_\alpha \rightarrow \Gamma$.
Given $\eta \in \KK^*_{af} \left( \coprod_\alpha B\Gamma^\prime_\alpha  \rightarrow B\Gamma \right)$ and a continuous map of pairs
$\nu :  
({\mathcal C}, {\mathcal C} - 
\Int(K)) \rightarrow (B\Gamma, B\Gamma^\prime)$,
we have
\begin{equation} \label{3.37}
\int_{\mathcal C} \widehat{A}(TM) \wedge 
\ch(\nu^* \eta) = 0.
\end{equation}
\end{corollary}

\begin{example} \label{3.38}
Let $\Gamma$ be a discrete group.
Put $\Gamma^\prime_1 = \Gamma^\prime_2 = \Gamma$. The exact sequence in $\KK$-theory gives
$\KK^0(\coprod_{\alpha} B\Gamma^\prime_\alpha \rightarrow B\Gamma) \cong
\KK^{-1}(B\Gamma)$.  Suppose that $\widehat{\eta} \in \KK^{-1}(B\Gamma)$ is almost flat
in the sense of Definition \ref{2.5}.  Then the corresponding relative element
$\eta \in \KK^0(\coprod_{\alpha} B\Gamma^\prime_\alpha \rightarrow B\Gamma)$ is almost flat,
which can be seen as follows.  Let ${\mathcal N}$, 
${\mathcal N}^\prime_1$ and ${\mathcal N}^\prime_2$ be compact manifolds-with-boundary, with smooth
maps $\mu_\alpha : {\mathcal N}^\prime_\alpha \rightarrow {\mathcal N}$. Let
$\nu : ({\mathcal N}, {\mathcal N}^\prime) \rightarrow (B\Gamma, B\Gamma^\prime)$ be a
continuous map.
Let $(V, \sigma)$ be as in Definitions \ref{2.4} and \ref{2.5}. We rename $\sigma$ as $\sigma_{odd}$.
Put $E = V \oplus V$, with the obvious $\Z_2$-grading. Define $\sigma_1 : \mu_1^* E^+
\rightarrow \mu_1^* E^-$ to be the identity on 
$\mu_1^* V$ and define $\sigma
: \mu_2^* E^+ \rightarrow \mu_2^* E^-$ to be 
$\mu_2^* \sigma_{odd}$.
Then $(E, \{\sigma_\alpha\})$ satisfies the requirements of Definition \ref{3.9}.
Multiplying by $\Z$, there is a similar construction of elements of
$\KK^{-1}_{af}(\coprod_{\alpha} B\Gamma^\prime_\alpha \rightarrow B\Gamma)$.

Now let $M$ be a complete Riemannian spin manifold with positive scalar curvature.
Let $Z$ be a compact separating hypersurface.
Given $\nu : M \rightarrow B\Gamma$, Theorem \ref{3.27}
implies that
$\int_{Z} \widehat{A}(TZ) \wedge \ch(\nu \big|_{Z}^* \widehat{\eta}) = 0$.
This can be compared with \cite[Theorem 1.5]{Cecchini-Rade-Zeidler (2023)}, 
\cite[Theorem 1.1]{Chen-Liu-Shi-Zhu (2021)} and
\cite[Corollary 6.8 and Theorem 6.12]{Gromov-Lawson (1983)}.
\end{example}

\begin{example} \label{3.39}
As an extension of the previous example,
let $\Gamma^\prime_1$, $\Gamma^\prime_2$ and $\Gamma_3$ be discrete groups, with injective homomorphisms
$\Gamma_3 \rightarrow \Gamma^\prime_1$ and $\Gamma_3 \rightarrow \Gamma^\prime_2$. Put
$\Gamma = \Gamma^\prime_1 \star_{\Gamma_3} \Gamma^\prime_2$, the amalgamated free product.  There are
injective homomorphisms $\Gamma^\prime_1 \rightarrow \Gamma$ and $\Gamma^\prime_2 \rightarrow \Gamma$.
The exact sequence in $\KK$-theory gives
$\KK^*(\coprod_{\alpha} B\Gamma^\prime_\alpha \rightarrow B\Gamma) \cong
\KK^{*-1}(B\Gamma_3)$. Given an almost flat element $\widehat{\eta} \in \KK^{*-1}(B\Gamma_3)$
in the sense of Definition \ref{2.5}, it is not immediate that the corresponding element
$\eta \in \KK^*(\coprod_{\alpha} B\Gamma^\prime_\alpha \rightarrow B\Gamma)$ is always almost flat,
but suppose that it is.
Then we can conclude the following.
Let $M$ be a complete Riemannian spin manifold with positive scalar curvature.
Let $Z$ be a compact separating connected $\pi_1$-injective hypersurface, dividing $M$ into pieces $M_1$ and $M_2$, with
$M_1 \cap M_2 = Z$.
Put $\Gamma_3 = \pi_1(Z, z_0)$, $\Gamma^\prime_1 = \pi_1(M_1, z_0)$ and 
$\Gamma^\prime_2 = \pi_1(M_2, z_0)$.
Let $\nu_Z : Z \rightarrow B\Gamma_3$ be the classifying map for the universal cover of $Z$.
Theorem \ref{3.27} implies that
$\int_{Z} \widehat{A}(TZ) \wedge \ch (\nu_Z^* \widehat{\eta}) = 0$.
\end{example}

\begin{corollary} \label{3.40}
Suppose that we are given
\begin{itemize}
\item A compact spin manifold-with-boundary $Y$, with boundary components $\{Y^\prime_\alpha\}$,
\item Discrete groups $\Gamma$ and $\{ \Gamma^\prime_\alpha \}$, with
homomorphisms $h_\alpha : \Gamma^\prime_\alpha \rightarrow \Gamma$, and
\item A continuous map $\nu: (Y, \partial Y) \rightarrow (B\Gamma, B\Gamma^\prime)$.
\end{itemize}
Suppose that the interior of $Y$ has a complete Riemannian metric
with positive scalar curvature.
Then for any $\eta \in \KK^*_{af}(\coprod_\alpha B\Gamma^\prime_\alpha \rightarrow B\Gamma)$,
we have $\int_Y \widehat{A}(TY) \wedge \ch(\nu^* \eta) = 0$.
\end{corollary}
\begin{proof}
Taking $K$ and ${\mathcal C}$ to be sufficiently large compact subsets of $\Int(Y)$ that
are diffeomorphic to $Y$, the
corollary follows from Theorem \ref{3.27}.
    \end{proof}

    \begin{corollary} \label{3.41}
    Let $Y$ be a compact connected spin manifold-with-boundary, with boundary components $\{Y^\prime_\alpha\}$.
    Put $\Gamma = \pi_1(Y)$ and $\Gamma^\prime_\alpha = \pi_1(Y_\alpha)$, where basepoints are handled as in 
    Section \ref{subsect3.2}.  Suppose that
    the interior of $Y$ has a complete Riemannian metric with positive scalar curvature,
    and
    $\KK^*_{af}(\coprod_\alpha B\Gamma^\prime_\alpha \rightarrow B\Gamma) \otimes \Q = 
    \KK^*(\coprod_\alpha B\Gamma^\prime_\alpha \rightarrow B\Gamma) \otimes \Q$.
    Let $\nu : (Y, \partial Y) \rightarrow (B\Gamma, B\Gamma^\prime)$ be the 
    classifying map.
    Then $\nu_*[Y, \partial Y]$ vanishes in $\HH_{\dim(Y)}(\coprod_\alpha B\Gamma^\prime_\alpha \rightarrow B\Gamma; \Q)$.
    \end{corollary}
    \begin{proof}
    Given Corollary \ref{3.40}, the proof is similar to that of Corollary \ref{2.7}.
    \end{proof}

\section{Obstructions from almost flat stable  relative bundles} \label{sect4}

In this section we give local obstructions to positive scalar curvature using almost flat stable relative bundles.
We use it to give obstructions for a manifold to admit a complete finite volume metric with positive scalar curvature.

By definition, a relative $\KK$-theory class has a better chance to be an almost flat stable class than to be
an almost flat class.  Example \ref{4.3} below shows that for pairs of aspherical manifolds, it is indeed possible that
a relative $\KK$-theory class fails to be an almost flat class, but is an almost flat stable class. Roughly speaking,
for pairs $(B\Gamma, B\Gamma^\prime)$, we can only expect relative $K$-theory classes to be
almost flat when the homomorphism $\Gamma^\prime \rightarrow \Gamma$ is injective, whereas
the notion of {\em stable} almost flatness has a wider range of validity.  

\begin{definition} \label{4.1}
Let ${\mathcal N}$ be a compact manifold-with-boundary.  Let ${{\mathcal N}^\prime}$
be a compact manifold-with-boundary with a smooth map
$\mu : {\mathcal N}^\prime \rightarrow {\mathcal N}$.
Put Riemannian metrics on ${\mathcal N}$ and
${{\mathcal N}^\prime }$.
Given $\beta \in \KK^0 \left( {\mathcal N}^\prime \rightarrow {\mathcal N}
\right)$, we say that $\beta$ is an almost flat stable class if for every
$\epsilon > 0$,
we can find
\begin{itemize}
    \item A $\Z_2$-graded Hermitian vector bundle $E$ on ${\mathcal N}$ and
    a Hermitian vector bundle $V$ on ${\mathcal N}^\prime$,
    \item 
    An isometric isomorphism $\sigma : \mu^* E^+  \oplus V
    \rightarrow \mu^* E^- \oplus V$ with
    $(E, {V}, {\sigma})$ representing $\beta$, and
    \item Hermitian connections $\nabla^E$ and $\nabla^{V}$ on $E$ and $V$, respectively,
    whose curvatures satisfies
    $\| F^E \| \le \epsilon$ and $\| F^{V} \| \le \epsilon$, and also 
    $\|\nabla \sigma \| \le \epsilon$.
\end{itemize}
\end{definition}

Let $\KK^0_{af,st} \left( {\mathcal N}^\prime \rightarrow 
{\mathcal N} \right)$ denote the almost flat stable classes.
The above definition of an almost flat stable class is essentially the same as that in
\cite[Definition 3.17]{Kubota (2022)}. The definition in \cite{Kubota (2022)} is effectively when 
$\mu$ is an embedding.

\begin{definition} \label{4.2}
Let $\Gamma$ be a discrete group and let $\{\Gamma^\prime_\alpha\}$ be a finite collection of
discrete groups, with homomorphisms $h_\alpha : \Gamma^\prime_\alpha \rightarrow \Gamma$.
Given $\eta \in \KK^{0} ( \coprod_\alpha B\Gamma^\prime_\alpha \rightarrow B\Gamma)$,
we say that $\eta$ is an almost flat stable class if 
for every compact manifold-with-boundary ${\mathcal N}$, 
compact manifolds-with-boundary $\{ {\mathcal N}^\prime_\alpha \}$, smooth maps
$\mu_\alpha : {\mathcal N}^\prime_\alpha \rightarrow {\mathcal N}^\prime$ and
continuous map of pairs
$\nu : ({\mathcal N}, {\mathcal N}^\prime) \rightarrow (B\Gamma, B\Gamma^\prime)$,
we have $\nu^* \eta \in \KK^0_{af,st} \left( \coprod_\alpha {\mathcal N}^\prime_\alpha \rightarrow 
{\mathcal N} \right)$.
\end{definition}

Let $\KK^{0}_{af,st}(\coprod_\alpha B\Gamma^\prime_\alpha \rightarrow B\Gamma)$ denote the almost flat stable classes.
Note that if $B\Gamma$ and $\{B\Gamma^\prime_\alpha\}$ are smooth manifolds then
$\eta \in \KK^{0}(\coprod_\alpha B\Gamma^\prime_\alpha \rightarrow B\Gamma)$ is an almost flat stable class in
the sense of Definition \ref{4.2} if and
only if it is an almost flat stable class in the sense of Definition \ref{4.1}.

Define the almost flat stable elements
$\KK^{-1}_{af,st}(\coprod_\alpha B\Gamma^\prime_\alpha \rightarrow B\Gamma)$ to be the elements of
$\KK^{-1}(\coprod_\alpha B\Gamma^\prime_\alpha \rightarrow B\Gamma)$
whose product with
$[B\Z] \in \KK^{-1}(B\Z)$ lies in 
$\KK^{0}_{af,st}(\coprod_\alpha B(\Gamma^\prime_\alpha \times \Z) \rightarrow 
B(\Gamma \times \Z))$.

It is conceivable that for all discrete groups $\Gamma$ and $\{\Gamma^\prime_\alpha\}$, and 
homomorphisms $h_\alpha : \Gamma^\prime_\alpha \rightarrow \Gamma$, one has
$\KK^{*}_{af,st}(\coprod_\alpha B\Gamma^\prime_\alpha \rightarrow B\Gamma) \otimes \Q = 
\KK^{*}(\coprod_\alpha B\Gamma^\prime_\alpha \rightarrow B\Gamma) \otimes \Q$, at least if $B\Gamma$ and the
$B \Gamma^\prime_\alpha$'s are finite CW-complexes.

\begin{example} \label{4.3}
Let $\Gamma$ be the trivial group and let
$\Gamma^\prime_1$ be any discrete group.  The exact sequence in $\KK$-theory gives
$\KK^*(B\Gamma^\prime_1 \rightarrow B\Gamma) \cong
\KK^{*-1}(B\Gamma^\prime_1)/\KK^{*-1}(\pt)$. Suppose that $\widehat{\eta} \in 
\KK^{-1}(B\Gamma^\prime_1)$ is almost flat in the sense of Definition \ref{2.5}.
Let ${\mathcal N}$ and ${\mathcal N}^\prime_1$ be compact Riemannian manifolds-with-boundary,
with a smooth map $\mu_1 : {\mathcal N}^\prime_1 \rightarrow {\mathcal N}$.
Let $\nu : ({\mathcal N}, {\mathcal N}^\prime_1) \rightarrow (B\Gamma,  B\Gamma^\prime_1)$
be a continuous map.
With reference to Definition \ref{2.4}, let $(V, \widehat{\sigma})$ be an almost flat pair that
represents $(\nu \big|_{{\mathcal N}^\prime_1})^* \widehat{\eta}$. With reference to
Definition \ref{4.1}, let $E^\pm$ be trivial vector bundles on ${\mathcal N}$ with a product
connection.  Put $\sigma = \Id \oplus \widehat{\sigma}$. Then $(E, V, \sigma)$ 
represents an almost flat stable class in $\KK^0({\mathcal N}^\prime_1 \rightarrow
{\mathcal N})$. Hence $\KK^0_{af,st}(B\Gamma^\prime_1 \rightarrow B\Gamma) = 
\KK^0(B\Gamma^\prime_1 \rightarrow B\Gamma)$.  There is a similar statement for
$\KK^{-1}_{af,st}(B\Gamma^\prime_1 \rightarrow B\Gamma)$. Compare with Example \ref{3.13}.
\end{example}

In this section we make the following assumption.
\begin{assumption} \label{4.4}
Given $r_0, D > 0$, put $r_0^\prime = \frac{1}{256} r_0^2 D^2$
and $D^\prime = D + \frac{32}{r_0 D}$. 
Let $M$ be a Riemannian spin 
manifold, possibly incomplete. Let $K$ be a compact codimension-zero submanifold-with-boundary in $M$.
Suppose that 
\begin{itemize}
\item The distance neighborhood $N_{D^\prime}(K)$ lies in a compact submanifold-with-boundary
${\mathcal C}$ of $M$,
\item $R > 0$ on $K$,
\item $R \ge r_0$ on $N_{D}(K) - K$,
\item $R \ge - r_0^\prime$ on $N_{D^\prime}(K) - N_{D}(K)$, and
\item $\partial K$ has nonnegative mean curvature as seen from
$M-K$.
\end{itemize}
\end{assumption}

\begin{theorem} \label{4.5}
Suppose that Assumption \ref{4.4} holds.
We can identify $\KK^* \left( ({\mathcal C} - \Int(K))
\rightarrow {\mathcal C} \right)$ with
$\KK^*\left({\mathcal C},  {\mathcal C} - \Int(K) \right)$.
Then given $\beta \in \KK^{*}_{af,st} \left( 
{\mathcal C},  {\mathcal C} - \Int(K) \right)$,
we have
\begin{equation} \label{4.6}
\int_{\mathcal C} \widehat{A}(TM) \wedge \ch(\beta) = 0.
\end{equation}
\end{theorem}
\begin{proof}
Suppose first that $M$ is even dimensional.  Let $\epsilon > 0$ be a
small parameter, which we will adjust.  Let ${\mathcal N}^\epsilon$ and
$f$ be as in the proof of Theorem \ref{2.9}.
As in Definition \ref{4.1}, let $(E, V, \sigma)$ be a triple that represents $\beta$,
where $E$ is a $\Z_2$-graded Hermitian vector bundle on ${\mathcal N^\epsilon}$,
$V$ is a Hermitian vector bundle on ${\mathcal N}^\epsilon - \Int(K)$
and $\sigma : E^+ \big|_{{\mathcal N}^\epsilon - \Int(K)} \oplus V \rightarrow
E^- \big|_{{\mathcal N}^\epsilon - \Int(K)} \oplus V$ is an isometric isomorphism.
Put $W = V \oplus V$, a $\Z_2$-graded vector bundle, and
 put $\nabla^W =
\nabla^V \oplus \nabla^V$. On the complement of $\Int(K)$, put
$Z = E \oplus W$, with Hermitian connection
$\nabla^{\pm} = \nabla^{E_\pm} \oplus \nabla^{W_\pm}$.

Define $\Pi^\prime$ as in the top of (\ref{2.18}), acting on sections of
$(S \otimes W) \Big|_{\partial K}$. Define  
$\Pi$ as in the bottom of (\ref{2.18}), acting on sections of 
$(S \otimes Z) \Big|_{\partial {\mathcal N}^\epsilon}$.
We will define a Dirac-type operator $D$ acting on elements
\begin{equation} \label{4.7}
(\psi, \psi^\prime) \in C^\infty({\mathcal N}^\epsilon; S \otimes E) \oplus 
C^\infty({\mathcal N}^\epsilon - \Int(K); S \otimes W)
\end{equation}
satisfying the boundary conditions 
$(\Pi^\prime - I) \psi^\prime \Big|_{\partial K} = 0$ and
$(\Pi - I) (\psi \oplus \psi^\prime) \Big|_{\partial {\mathcal N}^\epsilon} = 0$.
On $K$, we take $D$ to be the standard Dirac operator on 
$C^\infty(K; S \otimes E)$. On 
${\mathcal N}^\epsilon - \Int(K)$, where $D$ acts on sections of $S \otimes Z$,
we take $D$ to be
\begin{equation} \label{4.8}
D = 
\begin{pmatrix}
D^+ & \epsilon_S f \sigma^{-1} \\
\epsilon_S f \sigma & D^-
\end{pmatrix}.
\end{equation}
Here $D^\pm$ is the Dirac operator on sections of $S \otimes Z^\pm$ and $\epsilon_S$ is the $\Z_2$-grading
operator on spinors on $M$. Then on ${\mathcal N}^\epsilon - \Int(K)$,
\begin{align} \label{4.9}
D^2 = 
& \begin{pmatrix}
(D^+)^2 + f^2 & 0 \\
0 & (D^-)^2 + f^2
\end{pmatrix} + 
\sqrt{-1} \epsilon_S c(df)
\begin{pmatrix}
0 &  \sigma^{-1} \\
 \sigma & 0
\end{pmatrix} + \\
& \sqrt{-1} \epsilon_S f
\begin{pmatrix}
0 & c(\nabla^{+} \circ \sigma^{-1} - \sigma^{-1} \circ \nabla^{-}) \\
 c(\nabla^{-} \circ \sigma - \sigma \circ \nabla^{+}) & 0
\end{pmatrix}. \notag
\end{align}
Let ${\mathcal D}$ denote the ensuing self-adjoint operator with dense domain (\ref{4.7}).
As in Lemma \ref{2.19}, 
there is an $\epsilon_0 > 0$ so that if $\epsilon < \epsilon_0$, and in addition
$\| F^\pm \| \le \epsilon$ and $\| \nabla \sigma \| \le \epsilon$, then
the kernel of ${\mathcal D}$ vanishes. The condition that $\partial K$ has
nonnegative mean curvature, as seen from $M-K$, arises when considering the
boundary term from $\partial K$ on the spinor $\psi^\prime$
(defined on ${\mathcal N}^\epsilon - \Int(K)$), in the same
way that the condition on $\partial M$ arises in the proof of Theorem \ref{2.9}.

We compute the index of ${\mathcal D}$ along the lines of the
proofs of Lemma \ref{2.35} and Theorem \ref{3.27}.
We deform the Riemannian metric on ${\mathcal N}^\epsilon$ so that it is
a product near $\partial {\mathcal N}^\epsilon$, and a product in a neighborhood of
$\partial K$. Let us choose a small $\delta > 0$ so to parametrize a product neighborhood of
$\partial {\mathcal N}^\epsilon$ as $\partial {\mathcal N}^\epsilon \times [0,\delta]$,
and a product neighborhood of $\partial K$ in ${\mathcal N}^\epsilon - \Int(K)$ as
$\partial K \times [0, \delta]$.
We deform $f$ to zero. 
We deform $\sigma$ so that it pulls back from $\partial {\mathcal N}^\epsilon$ on
$\partial {\mathcal N}^\epsilon \times [0,\delta]$.
We deform $\nabla^V$ to be a product on $\partial K \times [0, \delta]$, with
$\nabla^W = \nabla^V \oplus \nabla^V$ there.
We deform $\nabla^Z$ so that it is a direct sum $\nabla^E \oplus \nabla^W$ 
on $\partial K \times [0, \delta]$.
We also deform $\nabla^Z$ to be a product on $\partial {\mathcal N}^\epsilon \times [0,\delta]$, with
$\nabla^{Z_-}  = \sigma \circ \nabla^{Z_+} \circ \sigma^{-1}$
there.

We can form a parametrix for the ensuing operator $D$ as follows. We take the double of $\partial K \times [0, \delta)$
to form $\partial K \times (- \delta, \delta)$ and we extend $W$ to the double, as $DW$. 
Let $DS$ denote the spinor bundle on the double.
We construct a 
$\Z_2$-invariant interior parametrix (for the Dirac operator) that sends sections of $DS \otimes DW$ with compact support in 
$\partial K \times (- \delta/2, \delta/2)$ to sections of $DS \otimes DW$ with compact support in
$\partial K \times (- \delta, \delta)$. Restricting the parametrix to 
$\Z_2$-invariant fields gives a parametrix that sends 
sections of $S \otimes W$  with compact support in 
$\partial K \times [0, \delta/2)$ to sections of $S \otimes W$ with compact support in
$\partial K \times [0, \delta)$, satisfying the boundary condition on $\partial K$.

Similarly, we use doubling to construct a parametrix for the Dirac operator that sends
sections of $S \otimes Z$ with compact support in  
$\partial {\mathcal N}^\epsilon \times [0, \delta/2)$
to sections of $S \otimes Z$ with compact support in  
$\partial {\mathcal N}^\epsilon \times [0, \delta)$, satisfying the boundary condition on 
$\partial {\mathcal N}^\epsilon$.
Next, we use an interior parametrix to
construct a parametrix for the Dirac operator that sends sections of $S \otimes E$ with compact support in
$K \cup (\partial K \times [0, \delta/2)$ to sections of $S \otimes E$ with compact support in
$K \cup (\partial K \times [0, \delta)$. Finally, we use an interior parametrix to construct a
parametrix for the Dirac operator that sends sections of $S \otimes Z$ with compact support in
${\mathcal N}^\epsilon - 
(K \cup (\partial K \times [0, \delta/8]) \cup (\partial {\mathcal N}^\epsilon \times [0, \delta/8]))$
to sections of $S \otimes Z$ with compact support in
${\mathcal N}^\epsilon - 
(K \cup (\partial K \times [0, \delta/4]) \cup (\partial {\mathcal N}^\epsilon \times [0, \delta/4]))$.
Using a partition of unity, we glue these parametrices together to form a
parametrix for $D$ on ${\mathcal N}^\epsilon$, satisfying the boundary conditions on
$\partial K$ and ${\mathcal N}^\epsilon$.

From the construction of the parametrix, the index of $D$
is 
\begin{align}
& \int_{\partial K \times [0, \delta/4)} \widehat{A}(T{\mathcal N}^\epsilon) \wedge \tr_s \left( e^{ \frac{i}{2\pi} (\nabla^W)^2} \right) +
\int_{K \cup (\partial K \times [0, \delta/4))} \widehat{A}(T{\mathcal N}^\epsilon) \wedge \tr_s \left( e^{ \frac{i}{2\pi} (\nabla^E)^2} \right) + \\
& \int_{{\mathcal N}^\epsilon - 
(K \cup (\partial K \times [0, \delta/4]))}
\widehat{A}(T{\mathcal N}^\epsilon) \wedge \tr_s \left( e^{ \frac{i}{2\pi} (\nabla^Z)^2} \right) = \notag \\
& \int_{K} \widehat{A}(T{\mathcal N}^\epsilon) \wedge \tr_s \left( e^{ \frac{i}{2\pi} (\nabla^E)^2} \right) + 
\int_{{\mathcal N}^\epsilon - K}
\widehat{A}(T{\mathcal N}^\epsilon) \wedge \tr_s \left( e^{ \frac{i}{2\pi} (\nabla^Z)^2} \right). \notag
\end{align}
As in Subsection \ref{subsect3.1}, let $G$ be a vector bundle on ${\mathcal N}^\epsilon - \Int(K)$ so that
$V \oplus G$ is isomorphic to a trivial bundle on ${\mathcal N}^\epsilon - \Int(K)$, which we extend to a
trivial bundle $T$ on ${\mathcal N}^\epsilon$.
Define a $\Z_2$-graded vector bundle $F$ on ${\mathcal N}^\epsilon$ by saying that
$F \big|_K = E|_K \oplus (T \oplus T)|_K$ and $F \big|_{{\mathcal N}^\epsilon - \Int(K)} =
Z \oplus (G \oplus G)$. Let $\nabla^G$ be an arbitrary connection on $G$ and let $\nabla^T$ be a connection on
$T$ that extends $\nabla^V \oplus \nabla^G$. Then
\begin{equation}
\int_{K} \widehat{A}(T{\mathcal N}^\epsilon) \wedge \tr_s \left( e^{ \frac{i}{2\pi} (\nabla^E)^2} \right) =
\int_{K} \widehat{A}(T{\mathcal N}^\epsilon) \wedge \tr_s \left( e^{ \frac{i}{2\pi} (\nabla^F)^2} \right)
\end{equation} 
because of the cancellation between the two $T$-factors, and
\begin{equation}
\int_{{\mathcal N}^\epsilon - K}
\widehat{A}(T{\mathcal N}^\epsilon) \wedge \tr_s \left( e^{ \frac{i}{2\pi} (\nabla^Z)^2} \right) =
\int_{{\mathcal N}^\epsilon - K}
\widehat{A}(T{\mathcal N}^\epsilon) \wedge \tr_s \left( e^{ \frac{i}{2\pi} (\nabla^F)^2} \right)
\end{equation} 
because of the cancellation between the two $G$-factors. Hence the index of $D$ equals
$\int_{{\mathcal N}^\epsilon} \widehat{A}(T{\mathcal N}^\epsilon) \wedge \tr_s \left( e^{ \frac{i}{2\pi} (\nabla^F)^2} \right)$. Because $\tr_s \left( e^{ \frac{i}{2\pi} (\nabla^F)^2} \right)$ vanishes in a neighborhood of
$\partial {\mathcal N}^\epsilon$,
we can extend it by zero to 
${\mathcal C}$ and conclude that the index of $D$ is $\int_{{\mathcal C}} \widehat{A}(TM) \wedge \tr_s \left( e^{ \frac{i}{2\pi} (\nabla^F)^2} \right)$. This equals
$\int_{\mathcal C} \widehat{A}(TM) \wedge \ch(\beta)$.  

If $M$ is odd dimensional, consider $M^\prime = M \times S^1$ 
and $K^\prime = K \times S^1$.  Then we reduce to the even dimensional case.
\end{proof}

\begin{corollary} \label{4.10}
Suppose that Assumption \ref{4.4} holds.
Let $\Gamma$ be a discrete group and let $\{\Gamma^\prime_\alpha\}$ be a finite collection of
discrete groups, with homomorphisms $h_\alpha : \Gamma^\prime_\alpha \rightarrow \Gamma$.
Given $\eta \in \KK^*_{af,st} \left( \coprod_\alpha B\Gamma^\prime_\alpha  \rightarrow B\Gamma \right)$ and a continuous map of pairs
$\nu :  
({\mathcal C}, {\mathcal C} - 
\Int(K)) \rightarrow (B\Gamma, B\Gamma^\prime)$,
we have
\begin{equation} \label{4.11}
\int_{\mathcal C} \widehat{A}(TM) \wedge 
\ch(\nu^* \eta) = 0.
\end{equation}
\end{corollary}

\begin{corollary} \label{4.12}
Suppose that Assumption \ref{4.4} holds, with ${\mathcal C}$ connected.
Let $\{({\mathcal C} - \Int(K))_\alpha\}$ be the connected components of ${\mathcal C} - \Int(K)$.
Put $\Gamma = \pi_1({\mathcal C})$ and $\Gamma^\prime_\alpha = \pi_1 \left( ({\mathcal C} - \Int(K))_\alpha \right)$,
where 
basepoints are handled as in Section \ref{subsect3.2}.
Let $\nu : ({\mathcal C}, {\mathcal C} - \Int(K)) \rightarrow (B\Gamma, B\Gamma^\prime)$ be the
classifying map.

If
$\KK^*_{af,st} \left( \coprod_\alpha B\Gamma^\prime_\alpha  \rightarrow B\Gamma \right) \otimes \Q =
\KK^* \left( \coprod_\alpha B\Gamma^\prime_\alpha  \rightarrow B\Gamma \right) \otimes \Q$
then $\nu_*[{\mathcal C}, {\mathcal C} - \Int(K)]$ vanishes in
$\HH_{\dim(M)} \left( \coprod_\alpha B\Gamma^\prime_\alpha  \rightarrow B\Gamma \right) \otimes \Q$.
\end{corollary}
\begin{proof}
Given Corollary \ref{4.10}, the proof is similar to that of Corollary \ref{2.7}.
\end{proof}

\begin{corollary} \label{4.13}
Let $M$ be a complete finite volume Riemannian spin manifold with positive scalar curvature and
 $\dim(M) \le 7$.
Then there is an exhaustion $K_1 \subset K_2 \subset \ldots$ of $M$ by compact submanifolds so that for each $j > i \ge 1$ and each $\beta \in \KK^*_{af,st} \left( K_j, K_j - \Int(K_i) \right)$, we have
$\int_{K_j} \widehat{A}(TM) \wedge \ch(\beta) = 0$.
\end{corollary}
\begin{proof}
Proposition \ref{2.41} gives a sequence $\{K_i\}$ of compact submanifolds that exhaust $M$ so that
$\partial K_i$ has nonnegative mean curvature as seen from $M - K_i$. Given $i$, if $j$ is
large enough then we can take ${\mathcal C} = K_j$ and satisfy Assumption \ref{4.4}.  
After passing to a subsequence of the $K_i$'s, the corollary follows from Theorem \ref{4.5}.
\end{proof}

\begin{corollary} \label{4.14}
Suppose that we are given
\begin{itemize}
\item A compact spin manifold-with-boundary $Y$, with boundary components $\{Y^\prime_\alpha\}$,
\item Discrete groups $\Gamma$ and $\{ \Gamma^\prime_\alpha \}$, with
homomorphisms $h_\alpha : \Gamma^\prime_\alpha \rightarrow \Gamma$, and
\item A continuous map $\nu: (Y, \partial Y) \rightarrow (B\Gamma, B\Gamma^\prime)$.
\end{itemize}
Suppose that $\dim(Y) \le 7$ and the interior of $Y$ has a complete finite volume Riemannian metric
with positive scalar curvature.
Then for any $\eta \in \KK^*_{af,st}(\coprod_\alpha B\Gamma^\prime_\alpha \rightarrow B\Gamma)$,
we have $\int_Y \widehat{A}(TY) \wedge \ch(\nu^* \eta) = 0$.
\end{corollary}
\begin{proof}
For small $\epsilon > 0$, let $T_\epsilon(\partial Y)$ be the $\epsilon$-neighborhood of
$\partial Y$ in $Y$. We can deform $\nu$ to a map $\widehat{\nu} : (Y, T_\epsilon(\partial Y)) 
\rightarrow (B\Gamma, B\Gamma^\prime)$.
Let $\{K_i\}$ be as in Corollary \ref{4.13}.  If $i$ is sufficiently large then for $j > i$, the class
$\ch \left( \widehat{\nu}\Big|_{(K_j, K_j - \Int(K_i))}^* \eta \right)$ can be represented by a differential form with
support in $K_i$ that extends by zero to a representative of $\ch \left( {\nu}^* \eta \right) \in
\HH^*_c(\Int(Y); \Q) \cong \HH^*(Y, \partial Y; \Q)$. It follows from Corollary \ref{4.13}
that $\int_Y \widehat{A}(TY) \wedge \ch(\nu^* \eta) =
\int_{K_j} \widehat{A}(TY) \wedge \ch(\widehat{\nu}^* \eta) = 0$.
The corollary follows.
    \end{proof}

    \begin{corollary} \label{4.15}
    Let $Y$ be a compact connected spin manifold-with-boundary, with boundary components $\{Y^\prime_\alpha\}$.
    Put $\Gamma = \pi_1(Y)$ and $\Gamma^\prime_\alpha = \pi_1(Y_\alpha)$, where basepoints are handled as in 
    Section \ref{subsect3.2}.  Suppose that
    \begin{itemize}
    \item $\Int(Y)$ has a complete finite volume metric with positive scalar curvature,
    \item $\dim(Y) \le 7$ and
    \item $\KK^*_{af,st}(\coprod_\alpha B\Gamma^\prime_\alpha \rightarrow B\Gamma) \otimes \Q = 
    \KK^*(\coprod_\alpha B\Gamma^\prime_\alpha \rightarrow B\Gamma) \otimes \Q$.
    \end{itemize}
    Then $\nu_*[Y, \partial Y]$ vanishes in $\HH_{\dim(Y)}(\coprod_\alpha B\Gamma^\prime_\alpha \rightarrow B\Gamma; \Q)$.
    \end{corollary}
    \begin{proof}
    Given Corollary \ref{4.14}, the proof is similar to that of Corollary \ref{2.7}.
    \end{proof}

\begin{example} \label{4.16}
Suppose that $\partial Y$ is connected.  If $\pi_1(Y) = \{e\}$ then there is some intersection with
Example \ref{2.46}.  If $\pi_1(\partial Y) = \{e\}$ then there is some intersection with 
Example \ref{3.34}, when the manifold $Z$ there is a copy of $\R^n$.
If $\pi_1(Y)$ and $\pi_1(\partial Y)$ are both nontrivial free abelian groups then Corollary \ref{4.15}
applies, which goes beyond the previous examples.
\end{example}

\appendix
\section{} \label{sectA}

In this appendix we discuss the relationship between simplicial volume and the index theoretic results
in the body of the paper.  In Subsection \ref{subsectA.1} we review simplicial volume and $l_1$-homology.
In Subsection \ref{subsectA.2} we consider the simplicial volume for closed Riemannian manifolds with almost
nonnegative scalar curvature, in two cases.  The first case is when there is an upper diameter bound and
a lower bound on the curvature operator.  The second case is when the volume is normalized to one and
there is a double sided curvature bound.  Our results in the second case are inconclusive.

\subsection{Review of simplicial volume} \label{subsectA.1}

Let $M$ be a closed compact connected oriented $n$-manifold.  The simplicial volume of $M$ 
\cite{Gromov (1982)} is
\begin{equation} \label{A.1}
\parallel M \parallel = \inf \left\{ \sum_i |a_i| \: : \: 
\left[ \sum_i a_i c_i \right] = [M] \right\},
\end{equation}
where $[M] \in \HH_n(M; \R)$ is the fundamental class in singular homology, each $c_i$ is a singular simplex
$c_i : \triangle^n \rightarrow M$ and the sum is finite. One knows that
$\parallel M \parallel$ is a homotopy invariant and just depends on the image 
$\nu_* [M] \in \HH_n(B\pi_1(M); \R)$ of $[M]$ under the
classifying map $\nu : M \rightarrow B\pi_1(M)$. 

In what follows we will refer to countable CW-complexes, although the statements will
apply to spaces that are homotopy equivalent to countable CW-complexes.
If $X$ is a countable CW-complex then the $l_1$-chains of $X$ are defined by
\begin{equation} \label{A.2}
C^{l_1}_k(X) = \left\{ \sum_i a_i c_i : \sum_i |a_i| < \infty \right\},
\end{equation}
where each $c_i$ is a singular $k$-simplex of $X$. The usual boundary operator gives a map
$\partial : C^{l_1}_k(X) \rightarrow C^{l_1}_{k-1}(X)$. The corresponding $l_1$-homology
groups are $\HH^{l_1}_*(X) = \Ker(\partial)/\Image(\partial)$. They acquire quotient
seminorms $\| \cdot \|_{l_1}$. 
There is a natural map $\HH_*(X; \Q) \rightarrow \HH^{l_1}_*(X)$.
If $M$ is a compact oriented manifold then
$\parallel M \parallel = \| [M] \|_{l^1}$. Despite the name, $l_1$-homology is not
a homology theory, in that the excision axiom is not satisfied,
which causes complications.

If $X$ is connected and $\nu : X \rightarrow B\pi_1(X, x_0)$ is the
canonical map, defined up to homotopy, then $\nu_* :
\HH^{l_1}_*(X) \rightarrow \HH^{l_1}_*(B\pi_1(X, x_0))$ is an isometric isomorphism \cite{Loeh (2008)}.
 
If $X$ is a countable CW-complex and $Y$ is a subcomplex of $X$ then
there is an
inclusion of $C^{l_1}_k(Y)$ into $C^{l_1}_k(X)$ and
we can define the
relative $l_1$-chains by $C^{l_1}_k(X,Y) = C^{l_1}_k(X)/C^{l_1}_k(Y)$. There is an induced
seminorm on $C^{l_1}_k(X,Y)$ and also an induced seminorm on the homology groups
$\HH^{l_1}_*(X,Y)$ of the corresponding chain complex. 

More generally, let $X$ be a countable CW-complex and let $Y$ be a countable CW-complexes equipped with a
cellular map $\mu : Y \rightarrow X$. 
We define the relative $l_1$-homology groups 
$\HH_*^{l_1}(Y \rightarrow X)$ by the
algebraic mapping cone construction.
That is, the relative chains are
\begin{equation} \label{A.3}
C_k^{l_1}(Y \rightarrow X) =
C_k^{l_1}(X) 
 \oplus C_{k-1}^{l_1}(Y)
\end{equation}
with boundary operator
$\partial(c_k, c_{k-1}) = (\partial c_k + \mu_* c_{k-1}, - 
\partial c_{k-1})$. 

If $Y$ is a subcomplex of $X$ then the chain map
$C_k^{l_1}(Y \rightarrow X) \rightarrow 
C_k^{l_1}(X, Y)$, given by
$(c_k, c_{k-1}) \rightarrow c_k \mod C^{l_1}_k(Y)$, induces an
isomorphism $\HH_*^{l_1}(Y \rightarrow X) \cong \HH^{l_1}_*(X,Y)$
which is in fact isometric \cite[Lemma 5.1]{BBFIPP (2014)} although we won't need this.

\begin{proposition} \label{A.4}
Let $X$ be a compact path connected topological space.
Let $Y$ be a closed subset of $X$, with path components
$\{Y_\alpha\}$.  Suppose that $(X,Y)$ is homotopy equivalent to a pair
of countable CW-complexes. Then 
\begin{equation} \label{A.5}
    \HH^{l_1}_*(X, Y) \cong \HH^{l_1}_* \left( \coprod_\alpha B\pi_1(Y_\alpha, y_\alpha) \rightarrow B\pi_1(X, x_0) \right),
    \end{equation}
    an isomorphism of
topological
vector spaces.
\end{proposition}
\begin{proof}
We have short exact sequences
\begin{equation} \label{A.6}
    \HH^{l_1}_*(Y) \rightarrow 
    \HH^{l_1}_*(X) \rightarrow
    \HH^{l_1}_*(Y \rightarrow X) \rightarrow
    \HH^{l_1}_{*-1}(Y) \rightarrow
    \HH^{l_1}_{*-1}(X) 
    \end{equation}
    and
    \begin{align} \label{A.7}
   & \HH^{l_1}_* \left( \coprod_\alpha B\pi_1(Y_\alpha, y_\alpha) \right)\rightarrow    
    \HH^{l_1}_*(B\pi_1(X,x_0)) \rightarrow \\
   & \HH^{l_1}_* \left( \coprod_\alpha B\pi_1(Y_\alpha, y_\alpha) 
    \rightarrow B\pi_1(X,x_0)) \right) \rightarrow 
 \HH^{l_1}_{*-1} \left( \coprod_\alpha B\pi_1(Y_\alpha, y_\alpha) \right)  
\rightarrow
\notag \\
    & \HH^{l_1}_{*-1}(B\pi_1(X,x_0)), \notag
\end{align}
along with arrows from the first sequence to the second sequence so that the
diagram commutes.  From \cite{Loeh (2008)}, the latter arrows (other than the middle one) are isomorphisms.
As $\HH^{l_1}_*(X, Y)$ is isomorphic to 
$\HH^{l_1}_*(Y \rightarrow X)$,
the proposition follows from the five lemma.
\end{proof}

\begin{remark} \label{A.8}
It would follow from the relative mapping theorem stated in \cite[Section 4.1]{Gromov (1982)}, along with \cite[Theorem 1.1]{Loeh (2008)}, that the isomorphism in Proposition \ref{A.4} is an isometry.
However, there doesn't seem to be a proof of the relative mapping theorem in the literature \cite[Remark 4.9]{Frigerio-Pagliantini (2012)}.
\end{remark}

\subsection{Simplicial volume and curvature} \label{subsectA.2}

\begin{conjecture} \label{A.9} \cite[Section 3.A]{Gromov (1986)}
For each $n \in \Z^+$, there is some $c_n > 0$ so that if $M$ is a 
compact connected oriented $n$-dimensional Riemannian manifold with
$R \ge - \sigma^2$ then $\parallel M \parallel \le c_n \sigma^n \vol(M)$.
\end{conjecture}

One can think of Conjecture \ref{A.9} two ways.  If we fix $\sigma = 1$ then 
$\parallel M \parallel$
would give an obstruction to volume-collapsing
with a lower scalar curvature bound.  If we fix $\vol(M) = 1$ then 
$\parallel M \parallel$
would give an obstruction for a manifold to have a 
Riemannian metric with
normalized volume and almost nonnegative scalar curvature.  We will think of it in the latter way.

The Ricci analog of Conjecture \ref{A.9} is known \cite{Gromov (1982)}. 
In addition there is a gap theorem; there is some $\epsilon = \epsilon(n) > 0$ so that if
$\vol(M) = 1$ and $\Ric \ge - \epsilon$ then $\parallel M \parallel = 0$. On the other hand,
Conjecture \ref{A.9} is not even known
when $\sigma = 0$. 

An analog of Conjecture \ref{A.9} for macroscopic scalar curvature, and a
gap result, are known \cite{Braun-Sauer (2021)}.

We look at whether Conjecture \ref{A.9} can be verified under some additional geometric bounds.
We first prove a gap result when there is an upper diameter bound and a lower bound on
the curvature operator.

\begin{proposition} \label{A.10}
Given $n \in \Z^+$ and $D, \Lambda < \infty$, there is an $\epsilon = \epsilon(n, D, \Lambda) > 0$ with the following property.
Let $M$ be a compact connected spin manifold of dimension $n$
whose fundamental group satisfies the Strong Novikov Conjecture for the maximal group $C^*$-algebra.  Suppose that $g$ is a Riemannian metric on $M$ so that 
$(M, g)$ has 
\begin{itemize}
\item Diameter bounded above by $D$, 
\item Curvature operator $\Rm$ bounded below by $- \Lambda$, and 
\item Scalar curvature $R$ bounded below by $-\epsilon$.
\end{itemize}
Then the simplicial volume $\parallel M \parallel$ vanishes. 
\end{proposition}
\begin{proof}
Suppose that the theorem is not true.  Then there is a sequence of
Riemannian $n$-manifolds $\{ (M_i, g_i) \}_{i=1}^\infty$ so that
\begin{itemize}
\item $\diam(M_i, g_i) \le D$, 
\item $\Rm(M_i, g_i)  \ge - \Lambda$, and 
\item $R(M_i, g_i) \ge - \frac{1}{i}$, but
\item $\parallel M_i \parallel > 0$.
\end{itemize}

The lower bound on $\Rm$ implies a lower bound on the Ricci curvature.
From \cite[Section 0.5]{Gromov (1982)}, the nonvanishing of the simplicial volume then implies that there is some $v_0 > 0$ so that $\vol(M_i, g_i) \ge v_0$
for all $i$. 

From \cite{Bamler-Cabezas-Rivas-Wilking (2019)}, there
are some $\tau > 0$ and $C < \infty$ so that there are Ricci flow solutions
$g_i(t)$, $t \in [0, \tau]$, with
\begin{itemize}
\item $g_i(0) = g_i$,
\item $\Rm(M_i, g_i(t)) \ge - C \Lambda$, and
\item $| \Rm(M_i, g_i(t)) | \le \frac{C}{t}$.
    \end{itemize}
    Put $g_i^\prime = g_i(\tau)$. The estimate $| \Rm(M_i, g_i(t)) | \le \frac{C}{t}$, along with Shi's derivative estimates, implies that we
    have uniform bounds on the $k^{th}$-covariant derivatives of
    $\Rm(M_i, g_i^\prime)$. It also implies a uniform upper bound
  $\diam(M_i, g_i^\prime) \le D^\prime$, by distance distortion estimates
  for the Ricci flow \cite[Section 27]{Kleiner-Lott (2008)}.
  By the monotonicity of scalar curvature under Ricci flow, we have
  $R(M_i, g_i^\prime) \ge - \frac{1}{i}$.

  Applying \cite[Section 0.5]{Gromov (1982)} again, there is some
  $v_0^\prime > 0$ so that 
  $\vol(M_i, g_i^\prime) \ge v_0^\prime$
for all $i$. After passing to a subsequence, there is a smooth limit
$\lim_{i\rightarrow \infty} (M_i, g_i^\prime) = (M_\infty, g_\infty)$,
where $M_\infty$ is diffeomorphic to each $M_i$, and
    $R(M_\infty, g_\infty) \ge 0$.

Either $(M_\infty, g_\infty)$ is Ricci flat or, after running the Ricci flow,
we can assume that it has positive scalar curvature.  If $(M_\infty, g_\infty)$ is Ricci flat then its fundamental group is virtually abelian and so
$\parallel M_\infty \parallel = 0$, which is a contradiction.  Hence
$(M_\infty, g_\infty)$ has positive scalar curvature. 

If $n$ is even then the Lichnerowicz
formula implies that the index of the Dirac operator on $(M_\infty, g_\infty)$ vanishes
in $\KK_0(C^*_{max} \pi_1(M_\infty))$.
Let $\nu : M_\infty \rightarrow B\pi_1(M_\infty)$ be the classifying map for the universal cover
of $M_\infty$; it is defined up to homotopy. Let $[M_\infty]_K \in \KK_n(M_\infty)$ denote the
fundamental class of $M_\infty$ in $K$-homology.
From the Strong Novikov Conjecture, the image $\nu_* [M_\infty]_K$ vanishes in $\KK_n(B\pi_1(M_\infty); \Q)$. Applying the Chern character gives
that $\nu_*(*\widehat{A}(TM_\infty))$ vanishes in $\HH_{*}(B\pi_1(M_\infty); \Q)$. As $* 1$ equals the fundamental class $[M_\infty] \in \HH_n(M_\infty; \Q)$, we
obtain that $\nu_*[M_\infty]$ vanishes in $\HH_n(B\pi_1(M_\infty); \Q)$. If $n$ is odd then
applying the preceding argument to $M_\infty \times S^1$ shows that 
$\nu_*[M_\infty]$ again vanishes in $\HH_n(B\pi_1(M_\infty); \Q)$.

Hence $\nu_*[M_\infty]$ vanishes
in the $l_1$-homology group $\HH_n^{l_1}(B\pi_1(M_\infty))$. From
\cite{Loeh (2008)}, the fundamental class $[M_\infty]$ now vanishes in
$\HH_n^{l_1}(M_\infty)$. Thus the simplicial volume $\parallel M_\infty
\parallel$ vanishes, which is a contradiction, since $M_\infty$ is diffeomorphic to $M_i$.
\end{proof}

\begin{corollary} \label{A.11}
Let $M$ be a compact connected spin manifold
whose fundamental group satisfies the Strong Novikov Conjecture for the maximal group $C^*$-algebra.
If $M$ admits a Riemannian metric with nonnegative scalar curvature then $\parallel M \parallel = 0$.
\end{corollary}

\begin{remark} \label{A.12}
Instead of using the Strong Novikov Conjecture, the conclusion of Proposition \ref{A.10} also holds
under the assumption that $\KK^*_{af}(B\pi_1(M)) \otimes \Q = \KK^*(B\pi_1(M)) \otimes \Q$.
This latter result also follows from Proposition \ref{A.10} when the Strong Novikov Conjecture for $C^*_{\max}$
is satisfied, using \cite{Hanke-Schick (2006)}.
\end{remark}

\begin{remark} \label{A.13}
One cannot remove the lower bound on the curvature operator in Proposition \ref{A.10}, as every closed
manifold of dimension greater than two admits Riemannian metrics for which $(\max |R|) \cdot \diam^2$ is
arbitrarily small \cite{Lohkamp (1999)}.
\end{remark}

We now look at what would be involved in proving a gap theorem saying
that for each $\Lambda < \infty$, there is some $\epsilon = \epsilon(n, \Lambda) > 0$ so that if a closed connected $n$-dimensional Riemannian manifold $(M,g)$ has $\vol(g) = 1$,
$R(g) \ge - \epsilon$ and $|\Rm|(g) \le \Lambda$ then $\parallel M \parallel = 0$.
We follow a contradiction argument along the lines of the proof of Proposition \ref{A.10}. The case when the
diameter stays bounded follows from Proposition \ref{A.10}, so we look at the case when the diameter goes
to infinity.  The next proposition says that there is a thick-thin decomposition in which the thick part
is modelled by a region in a complete noncompact finite volume Riemannian manifold with {\em positive}
scalar curvature.

\begin{proposition} \label{A.14}
Given $n \in \Z^+$, $v > 0$ and $\Lambda < \infty$
let
$\{(M_i,g_i)\}$ be a sequence of a closed connected $n$-dimensional Riemannian manifolds with 
\begin{itemize}
    \item $\vol(M_i, g_i) = 1$,
    \item $R(g_i) \ge - \frac{1}{i}$,
    \item $|\Riem|(g_i) \le \Lambda^2$ and
    \item $\sup_i \diam(M_i, g_i) = \infty$.
\end{itemize}
Put $M_i^{\ge v} = \{ m \in M_i \: : \: 
\vol(B(m, 1)) \ge v\}$.
Then after passing to a subsequence, for some $N$ there are
\begin{itemize}
\item A sequence of metrics $\{g_i^\prime\}$ that are 1.1-biLipschitz to 
$\{g_i\}$,
\item Points $\{m_{i,j}\}_{j=1}^N$ in $M_i^{\ge v}$ and
\item A collection $\{ (Z_j, \widehat{g}_j, p_j) \}_{j=1}^N$ of pointed connected complete finite volume noncompact
Riemannian $n$-manifolds with bounded
curvature and positive scalar curvature,
\end{itemize}
so that for each $1 \le j \le N$, we have $\lim_{i \rightarrow \infty} (M_i, g^\prime_i, m_{i,j}) = (Z_j, p_j)$ in the
pointed smooth topology, by approximants $\phi_{i,j} : (U_{i,j} \subset M_i) \rightarrow (V_{i,j} \subset Z_j)$
with $M_i^{\ge v} \subset \bigcup_{j=1}^N U_{i,j}$.
\end{proposition}
\begin{proof}
After passing to a subsequence, we can assume that $\lim_{i \rightarrow \infty} \diam(M_i, g_i) = \infty$.
From the double sided curvature bound, we can run the Ricci flow for a uniform time $\Delta > 0$
on each $M_i$, to obtain a metric $g_i^\prime$. Hence we can assume that
there are uniform bounds $|\nabla^I \Riem| \le C_{|I|}$ on the $(M_i, g_i^\prime)$'s.
Also, $R(g_i^\prime) \ge - \frac{1}{i}$.
By taking $\Delta$ small enough, we can ensure that each $g_i^\prime$ is $1.1$-biLipschitz to $g_i$.

For each $i$,
choose a maximal collection $B_i$ of disjoint unit $g_i$-balls
with center in $M_i^{\ge v}$. The
number of such balls is bounded above by $\lceil v^{-1} \rceil$.

If $B_i$ is empty for all large $i$ then we put
$N=0$ and stop.   If not then
after passing to a subsequence, for each $i$ we pick some $m_{i,1} \in M_i$
so that $B(m_{i,1}, 1)$ is an element of $B_i$.  After passing to a
subsequence, we can assume that there is a smooth pointed limit
$\lim_{i \rightarrow  \infty}(M_i, g^\prime_i, m_{i,1}) = 
(M_{\infty, 1}, g_{\infty, 1}, m_{\infty, 1})$; the fact that $\vol(B(m_{i,1}, 1) \ge v$ implies a
lower bound on the $g_i^\prime$-volume of the unit $g_i^\prime$-ball around $m_{i,1}$.

We now look whether after passing to a subsequence, for each $i$ we can find an element $B(m_{i,2}, 1)$ of $B_i$ so that
$\{d(m_{i,2}, m_{i,1})\}_{i=1}^\infty$ goes to infinity. If not, we stop and put $N=1$.
If so, we choose such new balls.  After passing to a subsequence, we
can assume that there is a smooth pointed limit
$\lim_{i \rightarrow  \infty}(M_i, g_i, m_{i,2}) = 
(M_{\infty, 2}, g_{\infty, 2}, m_{\infty, 2})$.

We now look whether after passing to a subsequence, for each $i$ we can find an element $B(m_{i,3}, 1)$ of $B_i$ so that
$\{d(m_{i,3}, \{m_{i,1}, m_{i,2}\})\}_{i=1}^\infty$ goes to infinity. If not, we stop and put $N=2$.
If so, we choose such new balls.  After passing to a subsequence, we
can assume that there is a smooth pointed limit
$\lim_{i \rightarrow  \infty}(M_i, g_i, m_{i,3}) = 
(M_{\infty, 3}, g_{\infty, 3}, m_{\infty, 3})$. 

We repeat the
process, which must terminate in at most $\lceil v^{-1} \rceil$ iterations. In the end,
we obtain a collection $\{(Z_j, \widehat{g}_j, p_j)\}_{j=1}^N$ of connected
pointed complete noncompact Riemannian manifolds with
\begin{itemize}
    \item $|\Riem(\widehat{g}_j)| \le C_0$ and
    \item $R(\widehat{g}_j) \ge 0$.
\end{itemize}
As $(Z_j, \widehat{g}_j)$ is the result of applying a Ricci flow
to a $W^{2,p}_{loc}$-regular Riemannian manifold,
the strong maximum principle implies that it is Ricci-flat or
has positive scalar curvature.  It has volume at most one.
Since a  noncompact complete finite volume
Riemannian manifold cannot be Ricci-flat \cite{Yau (1976)}, it follows that
$Z_j$ has positive scalar curvature. 

The conclusion of the proposition now holds; note that $M_i^{\ge v}$ is contained in the union of the
balls of radius $2$ with the
same centers as the elements of $B_i$.
\end{proof}

Proposition \ref{A.14}
implies that for large $i$, we can write $M_i = M_i^{thick} \cup M_i^{thin}$
where $M_i^{thick}$ and $ M_i^{thin}$ are $n$-dimensional submanifolds-with-boundary such that
\begin{itemize}
\item $\partial M_i^{thick} \subset \Int(M_i^{thin})$ and $\partial M_i^{thin} \subset \Int(M_i^{thick})$,
\item $M_i^{thick}$ is 1.1-biLipschitz equivalent to a union of regions in 
a finite collection of
pointed connected complete finite volume  noncompact
Riemannian $n$-manifolds with bounded
curvature and positive scalar curvature, and
\item Points in $M_i^{thin}$ are the centers of volume-collapsed unit balls.
\end{itemize}

Based on Corollary \ref{4.13}, we make the following hypothesis.

\begin{hypothesis} \label{A.15}
If $Z$ is a pointed connected complete finite volume oriented noncompact
Riemannian $n$-manifold with bounded
curvature and positive scalar curvature then there is an exhaustion $K_1 \subset K_2 \subset \ldots$ of $Z$ by
compact submanifolds with the following property.
Given $k \ge 1$, let $\{Y_{k,\alpha} \}$ be the connected components of
$K_{k+1} - \Int(K_k)$. Then
the pushforward $\nu_*[K_{k+1}, K_{k+1} - \Int(K_k)]$ vanishes in
$\HH_{\dim(Z)} \left( \coprod_\alpha B\pi_1(Y_{k,\alpha}) \rightarrow B\pi_1(K_{k+1}); \Q \right)$,
where basepoints are treated as in Section \ref{subsect3.2}.
\end{hypothesis}

Assuming Hypothesis \ref{A.15}, Proposition \ref{A.14} implies that we can write
$M_i = M_i^{thick} \cup M_i^{thin}$ where the fundamental class $[M_i^{thick},  
M_i^{thick} \cap M_i^{thin}]$ vanishes in $\HH^{l_1}_n(M_i^{thick},  
M_i^{thick} \cap M_i^{thin})$
and $M_i^{thin}$ is locally volume-collapsed.  From the lower curvature bound, there is a cover of
$M_i^{thin}$ by weakly amenable open sets with multiplicity at most $n$. From
\cite{Ivanov (2020), Loeh-Sauer (2020)}, the fundamental class $[M_i^{thin},
M_i^{thick} \cap M_i^{thin}]$ vanishes in $\HH^{l_1}_n(M_i^{thin},  
M_i^{thick} \cap M_i^{thin})$. We would now like to say that $\parallel M_i \parallel = 0$.
This amounts to a gluing property for simplicial volume. There are gluing results for the
simplicial volume but they do not seem to apply directly.

To explain the nature of the issue, given a fine triangulation of $M_i$, let $c$ be the sum of the (oriented)
$n$-simplices and let $c_{thick}$ be the sum of those in $\Int(M_i^{thick})$.
If the triangulation is sufficiently fine then
$\partial c_{thick}$ is a chain in $M_i^{thick} \cap M_i^{thin}$
and $(c_{thick}, - \partial c_{thick})$ represents the fundamental class in
$\HH^{l_1}_n(M_i^{thick},  
M_i^{thick} \cap M_i^{thin})$. Its vanishing means that there is an $l_1$-chain
$(d, e) \in C^{l_1}_{n+1}(M_i^{thick},  
M_i^{thick} \cap M_i^{thin})$ so that $\partial d + e = c_{thick}$. 
Then $c - \partial d = c - c_{thick} +e$ is $l_1$-homologous
to $c$ and has support in $M_i^{thin}$. There are vanishing results for
usual $n$-dimensional homology classes of $M_i^{thin}$ when mapped to $l_1$-homology
\cite{Ivanov (2020), Loeh-Sauer (2020)}
but they do not
seem to apply to $c - \partial d$.

To phrase the issue differently, we might have made a different choice $(d^\prime, e^\prime)$, which would change
$c - c_{thick} +e$ to $c - c_{thick} + e^\prime = c - c_{thick} +e + (e^\prime - e)$. As
$\partial e = \partial c_{thick}$, it follows that
$\partial(e^\prime - e) = 0$ and we can say that $e^\prime - e$ represents an element of
$\HH^{l_1}_n(M_i^{thick} \cap M_i^{thin})$ which lies in the image of the boundary map
$\HH^{l_1}_{n+1}(M_i^{thick}, M_i^{thick} \cap M_i^{thin}) \rightarrow \HH^{l_1}_n(M_i^{thick} \cap M_i^{thin})$.
It would seem necessary to be able to show that this element vanishes or, at least, has 
vanishing norm.  We know a bit more about the class of $e^\prime - e$
in $\HH^{l_1}_n(M_i^{thick} \cap M_i^{thin})$,
namely that under the isomorphism
$\HH^{l_1}_n(M_i^{thick} \cap M_i^{thin}) \cong \HH^{l_1}_n(B\pi_1(M_i^{thick} \cap M_i^{thin}))$
it arises from an ordinary homology class in $\HH_n(B\pi_1(M_i^{thick} \cap M_i^{thin}); \Q)$.


\begin{thebibliography}{99}

\bibitem{Adamson (1954)} I. Adamson,
`` Cohomology theory for non-normal subgroups and non-normal fields'', Proc. Glasgow Math.
Assoc. 2, p. 66-76 (1954)

\bibitem{Baer-Ballmann (2012)}
C. B\"ar and W. Ballmann, “Boundary value problems for elliptic differential operators of first order”,
in \underline{Surveys in Differential Geometry XVII}, International
Press, Boston, p. 1-78 (2012)

\bibitem{Baer-Hanke (2023)}
C. B\"ar and B. Hanke,
``Boundary conditions for scalar curvature'', in \underline{Perspectives in scalar curvature},
Vol. 2, 
World Scientific Publishing Co., Hackensack, NJ, 
p. 325-377 (2023)

\bibitem{Bamler-Cabezas-Rivas-Wilking (2019)} R. Bamler, E. Cabezas-Rivas and B. Wilking, "The Ricci flow under almost non-negative curvature conditions", Invent. Math. 217, p. 95-126 (2019)

\bibitem{Braun-Sauer (2021)} S. Braun and R. Sauer,
``Volume and macroscopic scalar curvature'',
Geom. Funct. Anal. 31, p. 1321-1376
(2021)

\bibitem{BBFIPP (2014)} M. Bucher, M. Burger, R. Frigerio, A. Iozzi, C. Pagliantini and  M. B. Pozzetti, ``Isometric embeddings in bounded
cohomology'', 
J. Topol. Anal. 6, p. 1-25 (2014)

\bibitem{Cecchini-Rade-Zeidler (2023)} S. Cecchini, D. R\"ade and R. Zeidler.
``Nonnegative scalar curvature on manifolds with at least two ends'', J. of Topology `6, p. 855-876 (2023)

\bibitem{Cecchini-Zeidler (2021)} S. Cecchini and R. Zeidler.
``The positive mass theorem and distance estimates in the spin setting'', 
Trans. Amer. Math. Soc. 377, p. 5271-5288  (2024)

\bibitem{Cecchini-Zeidler (2021b)} S. Cecchini and R. Zeidler.
``Scalar and mean curvature comparison via the Dirac operator'',
Geom. and Top. 28, p. 1167-1212 (2024)

\bibitem{Chang-Weinberger-Yu (2020)} S. Chang, S. Weinberger and G. Yu,
``Positive scalar curvature and a new index theory for
noncompact manifolds'', J. Geom. Phys. 149, 103575, 22 pp. (2020)

\bibitem{Chen-Liu-Shi-Zhu (2021)} J. Chen, P. Liu, Y. Shi and J. Zhu,
``Incompressible hypersurface, positive scalar curvature and positive mass theorem'',
preprint, https://arxiv.org/abs/2112.14442 (2021)

\bibitem{Chodosh-Li (2024)} O. Chodosh and C. Li,
``Generalized soap bubbles and the topology of manifolds with positive scalar curvature'',
Ann. Math. 199. p. 707-740 (2024)

\bibitem{Connes-Gromov-Moscovici (1990)} A. Connes, M. Gromov and H. Moscovici,
``Conjecture de Novikov et fibr\'es presque plats'', Comptes Rendus de l'Acad. des Sciences 310, p. 273-277 (1990)

\bibitem{Frigerio-Pagliantini (2012)} R. Frigerio and C. Pagliantini,
``Relative measure homology and continuous bounded cohomology of topological pair'',
Pac. J. Math. 257, p. 91-130 (2012)

\bibitem{Gromov (1982)} M. Gromov, ``Volume and bounded cohomology'',
Publ. Math. de l'IHES 56, p. 5-99(1982)

\bibitem{Gromov (1986)} M. Gromov, ``Large Riemannian manifolds'', in
\underline{Curvature and Topology of Riemannian Manifolds},
Lect. Notes in Math. 1201, Springer-Verlag, Berlin, p.108-122 (1986)

\bibitem{Gromov (1996)} M. Gromov,
``Positive curvature, macroscopic dimension, spectral gaps and higher signatures'', in
\underline{Functional analysis on the eve of the 21st century}, Volume II, Prog. Math. 132, p. 1-213, Birkh\"auser, Basel  (1996)


\bibitem{Gromov (2020)} M. Gromov, ``No metrics with positive scalar curvatures on aspherical 5-manifolds'',
preprint, https://arxiv.org/abs/2009.05332 (2020)

\bibitem{Gromov-Lawson (1980)}
M. Gromov and B. Lawson,
``The classification of simply connected
manifolds of positive scalar curvature'', Ann. Math 111, p. 423-434 (1980)

\bibitem{Gromov-Lawson (1983)} M. Gromov and B. Lawson,
``Positive scalar curvature and the Dirac operator on complete Riemannian manifolds'',
Publ. Math. IHES 58, p. 83-196 (1983)

\bibitem{Hanke-Schick (2006)} B. Hanke and T. Schick,
``Enlargeability and index theory'',
J. Diff. Geom. 74, p. 293-320
(2006)

\bibitem{Hunger (2019)} B. Hunger ``Almost flat bundles and homological
invariance of infinite K-area'', New York J. of Math. 25, p. 687-722 (2019)

\bibitem{Ivanov (2020)} N. Ivanov, ``Leray theorems for $l_1$-norms of infinite chains'',
preprint, https://arxiv.org/abs/2012.08690 (2020)

\bibitem{Karoubi (1978)} M. Karoubi, \underline{K-Theory}, Springer-Verlag, New York (1978)

\bibitem{Kleiner-Lott (2008)} B. Kleiner
and J. Lott, ``Notes on Perelman's papers'', Geom. Top. 12, p. 2587-2855 (2008)

\bibitem{Kubota (2022)} Y. Kubota,
``Almost flat relative vector bundles and the almost monodromy correspondence'',
J. Topol. Anal. 14, p. 353–382 (2022)

\bibitem{Loeh (2008)} C. L\"oh, ``Isomorphisms in $l^1$-homology'',
M\"unster J. of Math. 1, p. 237-266 (2008)

\bibitem{Loeh-Sauer (2020)} C. L\"oh and R. Sauer,
``Bounded cohomology of amenable covers via classifying spaces'',
Enseign. Math. 66, p. 151-172 (2020)

\bibitem{Lohkamp (1999)} J. Lohkamp.
``Scalar curvature and hammocks'', Math. Ann. 313, p. 385-407
(1999)

\bibitem{Lott (2021)} J. Lott,
``Index theory for scalar curvature on manifolds with boundary'', Proc. of the AMS 149, p. 4451-4459
(2021)

\bibitem{May (1999)} P. May, \underline{A concise course in algebraic topology},
Chicago Lectures in Mathematics, University of Chicago Press, Chicago (1999)

\bibitem{Mishchenko-Teleman (2005)}
A. Mishchenko and N. Teleman, ``Almost flat bundles and
almost flat structures'', Topol. Methods Nonlinear Anal. 26, p. 75-87 (2005)

\bibitem{Rosenberg (2007)} J. Rosenberg,
``Manifolds of positive scalar curvature: a progress report'',
in Surveys in Diff. Geom. 11,
\underline{Metric and comparison geometry}, International Press, Boston, p. 259-294 (2007)

\bibitem{Rosenberg-Weinberger} J. Rosenberg and S. Weinberger, 
``Positive scalar curvature on manifolds with boundary and their doubles'', Pure and
Appl. Math. Quart. 19, p.  2919-2950 (2023)

\bibitem{Schick-Seyedhosseini (2021)} T. Schick and M. Seyedhosseini,
``On an index theorem of Chang, Weinberger and Yu'', M\"unster J. Math. 14, p. 123-154 (2021)

\bibitem{Schoen-Yau (1979)} R. Schoen and S.T. Yau,
``On the structure of manifolds with positive scalar
curvature'', Manuscripta Math. 28, p. 159-183
(1979)

\bibitem{Segal (1968)} G. Segal, ``Classifying spaces and spectral sequences'',
Publ. Math. IHES, p. 105-112 (1968)

\bibitem{Skandalis (1991)} G. Skandalis,
``Approche de la conjecture de Novikov par la cohomologie cyclique (d'apr\`es A. Connes, M. Gromov et H. Moscovici)'',
S\'eminaire Bourbaki,
Ast\'erisque 201-203, p. 299-320 (1991)

\bibitem{Song (2023)} A. Song, ``A dichotomy for minimal hypersurfaces in manifolds thick at infinity'', 
Ann. Sci. \'Ec. Norm. Sup. 56, p. 1085-1134 (2023)

\bibitem{Takasu (1959)} S. Takasu, ``Relative homology and relative cohomology theory of groups'',
J. Fac. Sci. Univ. Tokyo. 
8, p. 75-110 (1959)

\bibitem{Wang-Zhang (2022)} X. Wang and W. Zhang,
``On the generalized Geroch conjecture for complete spin manifolds'',
Chinese Ann. Math. B43, p. 1143-1146 (2022)

\bibitem{Yau (1976)} S.-T. Yau, ``Some function-theoretic properties of complete Riemannian manifold and their applications to geometry''.
Indiana Univ. Math. J. 25, p. 659-670 (1976)

\end{thebibliography}
\end{document}